\documentclass[atmp,nocrop]{ipart_v1}
\renewcommand{\copyrightnotice}{}
\renewcommand{\clearemptydoublepage}{\clearpage}
\usepackage[T1]{fontenc}
\usepackage{lmodern}
\let\toprule\relax
\let\midrule\relax
\usepackage{amsmath,amssymb,amsthm,mathtools,booktabs,array,mathrsfs,needspace}
\usepackage{tikz-cd}
\usetikzlibrary{arrows.meta,positioning}
\hypersetup{colorlinks=true,linkcolor=black,citecolor=black,urlcolor=black}
\usepackage[nameinlink,noabbrev]{cleveref}
\allowdisplaybreaks[2]
\numberwithin{equation}{section}
\newtheorem{theorem}{Theorem}[section]
\newtheorem{proposition}[theorem]{Proposition}
\newtheorem{lemma}[theorem]{Lemma}
\newtheorem{corollary}[theorem]{Corollary}
\theoremstyle{definition}
\theoremstyle{remark}\newtheorem{remark}[theorem]{Remark}
\newcommand{\C}{\mathbb C}\newcommand{\Z}{\mathbb Z}\newcommand{\Q}{\mathbb Q}
\newcommand{\cO}{\mathcal O}\newcommand{\cN}{\mathcal N}\newcommand{\cB}{\mathcal B}
\newcommand{\cM}{\mathcal M}\newcommand{\cL}{\mathcal L}\newcommand{\cD}{\mathfrak D}

\newcommand{\cC}{\mathcal C}
\newcommand{\cK}{\mathcal K}\newcommand{\cR}{\mathcal R}\newcommand{\Gr}{\operatorname{Gr}}
\newcommand{\U}{\mathsf U}\newcommand{\T}{\mathsf T}\newcommand{\im}{\operatorname{im}}
\newcommand{\GL}{\operatorname{GL}}\newcommand{\SL}{\operatorname{SL}}
\newcommand{\Spec}{\operatorname{Spec}}\newcommand{\Sym}{\operatorname{Sym}}
\newcommand{\Fr}{\operatorname{Fr}}\newcommand{\Hom}{\operatorname{Hom}}
\newcommand{\Ext}{\operatorname{Ext}}
\newcommand{\wt}{\operatorname{wt}}\newcommand{\rk}{\operatorname{rk}}
\newcommand{\ch}{\operatorname{ch}}\newcommand{\Sh}{\operatorname{Sh}}
\newcommand{\pr}{\operatorname{pr}}
\newcommand{\eps}{\varepsilon}

\hypersetup{pdftitle={Uniform localization and boundary corrections for two-node quiver Coulomb branches},pdfauthor={Xiaomin Tang},pdfsubject={Coherent Coulomb branches and normalized difference operators}}
\title[Uniform localization and boundary corrections]{Uniform localization and boundary corrections for two-node quiver Coulomb branches}
\author[Xiaomin Tang]{Xiaomin Tang}
\begin{document}
\begin{abstract}
We construct simple coherent classes and compute their normalized difference
operators for the two-node quiver Coulomb branch with gauge group
$\GL_1\times\GL_2$ and $m-1$ framing columns at the second node.
The first magnetic charge is fixed at zero. For every $m\geq1$, every
support $(n,-1)$ with $n\geq1$, and every integral electric twist, we prove
that termwise saturation of the matter Koszul complex gives the simple
intermediate extension. A graded normalization formula expresses its
localization on every stratum as a finite binomial shuffle sum. The resulting
operators satisfy bar symmetry, positive bubbling order, paired root-residue
cancellation, and integral Laurent polynomial action. For $m=4$, we construct
the derived boundary correction on $(2,-2)$, retaining the nonsemisimple
odd-electric boundary object, and compute all electric twists. Finally,
matter Euler divisibility proves triangular uniqueness on arbitrary fixed
rank-two magnetic supports inside the actual integral Coulomb image.
A pole-free counterexample shows why residue conditions alone do not give
this uniqueness. Existence on general deeper mixed-sign supports remains open.
\end{abstract}
\maketitle
\noindent\textbf{2020 Mathematics Subject Classification.}
Primary 14M15, 14F08; Secondary 17B37, 81T60.\par
\noindent\textbf{Keywords.}
K-theoretic Coulomb branch, Koszul-perverse coherent sheaf, affine Grassmannian,
normalization cover, monopole bubbling, local cohomology.

\section{Introduction}
Quantized K-theoretic Coulomb branches admit both a coherent geometric description
and a realization by rational difference operators. The coherent description
retains the incidence equations, electric lines, and derived module structure;
abelianization makes the magnetic support and bubbling coefficients visible.
On a singular support, comparing these descriptions requires the localization of
the coherent module itself, not its ordinary fiber. Our purpose is to carry out
this comparison uniformly on a fixed infinite magnetic family, and then to
identify and remove the first boundary obstruction outside that family. For the
broader setting of Coulomb branches and symplectic duality, see the survey of
Kamnitzer~\cite{Kamnitzer}.

We address two related questions: when does saturation produce a simple
intermediate extension, and which algebraic conditions characterize its
normalized localization? The first-strip family $(n,-1)$ admits a uniform
answer for all $n\geq1$, all numbers of columns, and all integral electric
twists (Theorems~\ref{thm:main-geometry} and~\ref{thm:main}). For four columns,
the support $(2,-2)$ exhibits the first deeper mixed-sign obstruction:
a saturated coefficient fails to be maximal Cohen--Macaulay, and a derived
boundary subtraction is necessary (Theorem~\ref{thm:main-boundary}). A further
argument proves triangular uniqueness on every fixed rank-two support within
the actual integral Coulomb image (Theorem~\ref{thm:main-uniqueness}).
Throughout, ``full support'' means every stratum of one specified orbit
closure. Our existence results concern the stated families; the general
uniqueness theorem does not construct a complete canonical basis.

We work with the two-node datum
\begin{equation}\label{eq:datum}
 \begin{gathered}
 G=\GL(\C)\times\GL(V_0),\qquad V_0=\C^2,\\
 \mathbf N_m=\Hom(\C,V_0)\oplus\Hom(\C^{m-1},V_0),\quad m\geq1.
 \end{gathered}
\end{equation}
The first-node magnetic charge is always zero. The local matter representation
therefore has $m$ standard target columns, with source characters
$x_1,\ldots,x_m$; $x_1$ remains a gauge character, not a flavor specialization.
We write $x=\prod_{v=1}^m x_v$ and $e_r=e_r(x_1,\ldots,x_m)$, with $e_0=1$.
The target gauge rank is always two: the integer $n$ below is a magnetic charge,
not a variable rank.
\begin{figure}[ht]
\centering
\begin{tikzpicture}[>=Latex,node distance=24mm]
\node[draw,circle,minimum size=9mm] (g1) {$1$};
\node[draw,circle,minimum size=9mm,right=of g1] (g2) {$2$};
\node[draw,rectangle,minimum width=13mm,minimum height=9mm,right=of g2] (f) {$m-1$};
\draw[->] (g1)--(g2);\draw[->] (f)--(g2);
\node[below=3mm of g1,font=\small] {gauge};
\node[below=3mm of g2,font=\small] {gauge};
\node[below=3mm of f,font=\small] {framing};
\end{tikzpicture}
\caption{The fixed gauge dimension vector is $(1,2)$. Only the number of framing
columns varies; the first-node magnetic charge is zero.}
\label{fig:quiver}
\end{figure}
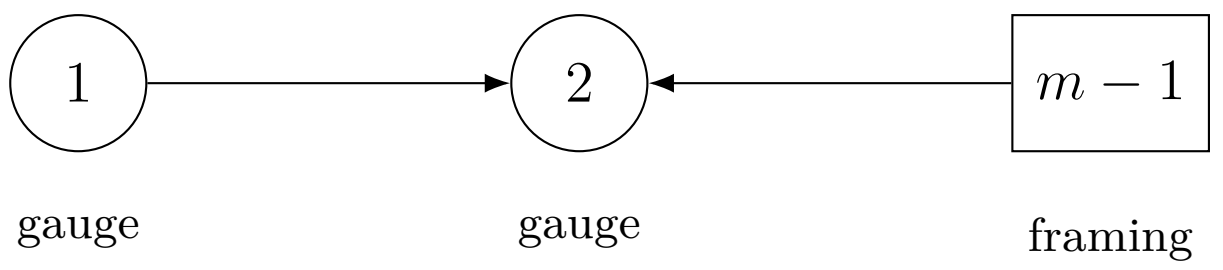

The four-column case is distinguished by its balanced $A_2$ quiver:
the gauge dimension vector is $(1,2)$ and the framing vector is $(0,3)$.
Indeed $0+2=2\cdot1$ and $3+1=2\cdot2$; the target excess for general $m$ is
$m-4$. This provides a natural setting for type-$A_2$ comparisons, but balance is
not a hypothesis of the first-strip theorem. The general Coulomb-branch
construction and its relation with affine Grassmannian slices are due to
Braverman--Finkelberg--Nakajima~\cite{BFN18,BFN19}. Cautis--Williams construct
the Koszul-perverse category and its canonical basis of simple classes~\cite{CW23},
extending the coherent Satake perspective developed in~\cite{CW19}.

For dimension vector $(1,2)$ and framing $(0,3)$, the reduced type-$A$
construction of Finkelberg--Tsymbaliuk is identified with an extended adjoint
integral form of quantum $\mathfrak{sl}_3$~\cite[\S4.3]{FT}. That comparison
uses its stated flavor reduction, covering groups, extended Cartan variables,
and complexified integral form; it does not identify our precise integral
character ring or simple-class labels. The differential-operator realization
of Gannon--Williams~\cite{GW} concerns an enlarged gauge datum; we do not
identify the present fixed $\GL_1\times\GL_2$ model with that enlarged
base-affine-space model.

The quiver-with-involution construction of Shen--Su--Wang--Xiong~\cite{SSWX}
also contains the datum~\eqref{eq:datum}: doubling the ordinary $A_2$ arrow
and placing the framing on its opposite copy recovers exactly $\mathbf N_m$,
without exterior-square or symmetric-square matter. This is a specialization
of their general construction, not their connected type-$A\mathrm{III}$ example.
Their quasi-minuscule localization formula~\cite[Theorem~2.2]{SSWX}, with its
homological precursor in~\cite[Theorem~2.4]{SSX}, gives a useful comparison at
magnetic charge $(0;1,-1)$. It concerns a resolution structure-sheaf pushforward,
up to an equivariant coefficient-ring term, rather than a Koszul-perverse simple
class. Thus it does not replace the intermediate-extension calculation below,
nor does it settle the higher first strip or the genuine $(2,-2)$ correction.
We make the datum and parameter comparison explicit in Section~\ref{sec:prelim}.

There is also a precise cluster-algebra setting for these computations.
Schrader--Shapiro~\cite[Theorem~5.13]{SS} identify the quantized $K$-theoretic
Coulomb algebra of a framed quiver without self-loops with its quantized
universally Laurent algebra, or upper cluster $\mathcal X$-algebra. Their framing
convention includes~\eqref{eq:datum}. Their separate cluster $\mathcal A$
construction in~\cite[Section~7]{SS}, however, assumes an unframed quiver with
vanishing first graph homology, and cannot be applied to our framed model on
that basis alone. Monoidal cluster categorifications are already known both
in coherent Satake~\cite{CW19} and with nonzero matter: Cautis--Williams construct a quantum monoidal seed and categorical exchange
structure of affine type $A_2^{(1)}$ for the centrally regraded
Koszul-perverse category of $\GL_2$ with standard matter $\C^2$
\cite[Theorems~1.5 and~8.10]{CW23}. Liu~\cite[Theorem~8.1]{LiuCluster} proves
such a categorification for certain star-shaped quiver theories, including
all of his central-rank-two cases. His datum has one framing column and
$r-1$ gauged one-dimensional legs. Our $m=2$ datum agrees with his
$(n,r)=(2,2)$ case, after identifying the scalar action conventions; our
$m=4$ datum, with one gauged leg and three framing columns, does not agree
with his $(2,4)$ case. Replacing gauge legs by framing is not an identification
of Coulomb algebras. These results supply concrete comparison frameworks,
not an automatic identification of the simple classes computed here with
cluster monomials or a theta basis.

Gaiotto's elementwise conditions for $K_{\mathfrak q}$-algebras~\cite{Gaiotto}
motivate the bar, positive-order, and residue tests below. The accompanying
KAlgebra software supplies a computational approach to canonical elements in
abelianized presentations~\cite{KAlgebra}. Our contribution is the construction
of matter-dependent geometric objects and their normalized localizations, not
a new general canonicalization algorithm. Satisfying the stated elementwise
tests is not a basis-wide comparison or a multiplication theorem.

In coherent Satake, Finkelberg--Fujita identify coherent IC classes with a
dual canonical basis~\cite{FF}; Liang constructs a neutral Tannakian
subcategory for adjoint matter in type $A$~\cite{Liang}.
Cautis--Vasserot~\cite{CV} relate affine quiver--Hecke simples to
Koszul-perverse simples for unframed finite simply laced quivers under a
strict-dominance hypothesis. These results have different matter or
categorical hypotheses from the framed datum here. In particular, balance of
$(v,w)=((1,2),(0,3))$ is not the strict-dominance hypothesis of~\cite{CV}.
We therefore do not infer a quantum-group basis identification or positivity
of normalized localizations from these comparisons.

For the geometry of the underlying varieties, Muthiah--Weekes describe
symplectic leaves of generalized affine Grassmannian slices~\cite{MW} and
embeddings of Kac--Moody slices through fundamental monopole operators~\cite{MWFMO}.
Webster gives local descriptions of Coulomb branches in the rational,
$K$-theoretic, and elliptic settings~\cite{Webster}. The calculations below
retain additional matter-incidence, electric-line, and dg information.
For an explicit $K$-theoretic comparison with different gauge data, see the
once-punctured-torus dualities of Allegretti--Shan~\cite{AS}.

We now state the four main results precisely. The first two concern the first
strip; the third identifies the first obstruction beyond it. The fourth
extends triangular uniqueness within the actual integral Coulomb image.
Fix $n\geq1$ and set $N=n+1$. The reduced support and its dominant strata are
\begin{equation}\label{eq:support}
 \begin{gathered}
 S=\overline{\cO_{(N-1,-1)}}_{\mathrm{red}},\qquad
 U=\cO_{(N-1,-1)},\\
 S_j=\cO_{\lambda_j},\quad \lambda_j=(N-1-j,j-1),
 \end{gathered}
\end{equation}
where $0\leq j\leq\lfloor N/2\rfloor$ and $\dim S_j=N-2j$.
Let $\iota:U\hookrightarrow S$ be the open inclusion and let $\ell$ be the
leading-image electric line. Its chosen upper fixed-point character is $b$.
All powers $k\in\Z$ are retained as integers; only a deck-character index will
be reduced modulo $N$.

\begin{theorem}\label{thm:main-geometry}
For every $m\geq1$, $N\geq2$, and $k\in\Z$, the intrinsic open matter Koszul
complex twisted by $\ell^k/x$ has a bounded coherent termwise saturated extension
$\cM_k$ on all of $S$. Every coefficient is maximal Cohen--Macaulay. The complex
carries the full ambient dg action and compatible algebraic positive-loop
equivariance. With the Cautis--Williams shift $[N/2-k]$, it is the simple
intermediate extension of its open incidence line, compatibly with the finite-stage
transitions defining the coherent Coulomb category.
\end{theorem}

The extension is constructed before localization. Its proof uses the actual cyclic
cover attached to $\ell^N$, not a root hypersurface chosen only from a fixed-point
weight. A rank-$(N-2)$ splitting reduces each exterior coefficient to finitely many
electric powers. Their maximal Cohen--Macaulay property follows from a finite
quasi-\'etale cover of a klt support. The differential, the ambient action, and the
algebraic jet-group linearization then extend by reflexive Hartogs arguments.
The strict boundary estimates are verified at every stratum, not only at the
first degeneration.

Let $H_i^{N,j}(t,\rho)$ denote the deck-character-$i$ Hilbert character of the
normalization of the two-row nilpotent slice of type $(N-j,j)$ in the universal
regular-nilpotent cover. The precise deck convention and the finite formula are
given in Section~\ref{sec:characters}. Let $u$ rotate $z$ with weight $u^2$ and
dilate matter with weight $u$. With the ordinary-to-symmetric normalization of
Section~\ref{sec:prelim}, put $\delta_j^N=a^jb^{N-j}$.

\begin{theorem}\label{thm:main}
The normalized localization of Theorem~\ref{thm:main-geometry} has symmetric
coefficients
\begin{equation}\label{eq:mainformula}
 \begin{aligned}
 f^{(k)}_{N-1,-1}(a,b)&=b^{k-m},\\
 f^{(k)}_{\lambda_j}(a,b)&=
 \frac{(-1)^ju^j}{x}\sum_{r=0}^m e_r\delta_j^{k-r}
 H_{r-k}^{N,j}\bigl(u^2,u^{N-2j}b/a\bigr),\quad 1\leq j\leq\lfloor N/2\rfloor.
 \end{aligned}
\end{equation}
Each $H_i^{N,j}$ is the $\binom Nj$-term descending-shuffle expression
\eqref{eq:shuffle}, including the shift $(N-\gcd(i,N))/2$ for nonzero sectors.
All fractional target characters cancel. Weyl reflection supplies the other
labels, and a central label is counted once.

The $N+1=n+2$ ordinary terms form the operator
\begin{equation}\label{eq:mainoperator}
 \T_{N,k}=\sum_{R=-1}^{N-1}B_R(a,b)\U_{R,N-2-R}.
\end{equation}
The symmetric coefficients are Weyl covariant and fixed by $u\mapsto u^{-1}$,
with all other characters fixed; the coefficient in stratum $j$ is $O(u^j)$.
The ordinary coefficients have only simple even-root poles and satisfy
\begin{equation}\label{eq:mainresidue}
 \rho_h(B_R)+\rho_h(B_{N-2-R+h})=0,\qquad
 \rho_h(B)=\left.(1-u^{-2h}b/a)B\right|_{b=u^{2h}a},
\end{equation}
where coefficients outside the support are zero. Consequently $\T_{N,k}$ preserves
\[
 \Z[u^{\pm1},x_1^{\pm1},\ldots,x_m^{\pm1}]
 [a^{\pm1},b^{\pm1}]^{\mathfrak S_2}.
\]
Within the integral even-root, simple-pole residue class on this support, it is the
unique bar-invariant positive-order triangular lift of its prescribed top
coefficient. After extraction of $(-1)^ju^j/x$, the common parameter denominator
in~\eqref{eq:positive} gives a nonnegative integral numerator, palindromic of
reciprocal $u$-degree $2j(N-j)$.
\end{theorem}

Graham's canonical-sheaf pushforward and vanishing compute the normalization
summands~\cite{Graham}; our graded specialization gives their finite shuffle
characters, including the nonprimitive-sector correction. After verifying the
categorical hypotheses and identifying the genuine coherent class, we apply
the Cautis--Williams framework and the Schrader--Shapiro image theorem~\cite{SS}.
This gives the normalized coefficients uniformly in $N$, $m$, and $k$.
At $k=m$ the leading coefficient is one.

For four columns the same electric-cover calculation also treats the same-sign
magnetic cones after recalculating their constant-rank matter incidence. It does
not extend unchanged to deeper mixed-sign charges. The first such support already
exhibits the difference between reflexive extension and intermediate extension.

For the next statement, let $\mathcal M_k$ be the termwise saturated
intrinsic matter complex on $\overline{\mathcal O}_{(2,-2)}$ constructed from
the electric/source twist $\ell^k/x$. Let $T_{k,l}$ be its exterior-degree-$l$
coefficient on the four-dimensional transverse slice, and let $0$ denote
the central point. Write $f^{\mathrm{sat},k}_{0,0}$
for the central symmetric coefficient of
$u^{2-2k}[\mathcal M_k[2-k]]/x$. The notation $\overline F$ means
$F(u^{-1},a,b,x_1,\ldots,x_4)$.

\Needspace{5\baselineskip}
\begin{theorem}\label{thm:main-boundary}
Let $m=4$ and $k\in\Z$. On the support
$(2,-2)>(1,-1)>(0,0)$, the following assertions hold.
\begin{enumerate}
\item The saturated complex satisfies
\[
 \operatorname{depth}_{0}T_{k,4}=2,
 \qquad \mathcal M_k[2-k]\in\operatorname{KP},
\]
where $\operatorname{KP}$ is the Koszul-perverse heart. There is a canonical
exact sequence in this heart
\begin{equation}\label{eq:main-boundary-sequence}
 0\longrightarrow\mathcal N_k\longrightarrow\mathcal M_k[2-k]
 \longrightarrow\operatorname{IC}_k\longrightarrow0.
\end{equation}
Here $\mathcal N_k$ is the maximal central subobject, formed from the full
row $H^2_{\{0\}}(T_{k,l})$ with its induced dg structure, and
$\operatorname{IC}_k$ is the simple intermediate extension of the open
incidence line. The object $\mathcal N_k$ is simple for even $k$ and
in general nonsemisimple for odd $k$.
\item The normalized class
\[
 \mathsf K_k=\frac{u^{2-2k}}{x}[\operatorname{IC}_k]
\]
has top symmetric coefficient $f^{(k)}_{2,-2}=b^{k-8}$, and its central
coefficient is
\begin{equation}\label{eq:main-boundary-correction}
 \begin{gathered}
 f^{(k)}_{0,0}=f^{\mathrm{sat},k}_{0,0}-\Delta_k,\\
 \Delta_k=\begin{cases}
 \displaystyle\frac{u^2(ab)^{(k-4)/2}}{x},&k\text{ even},\\[4pt]
 \displaystyle\frac{u(ab)^{(k-5)/2}}{x^2}
 \bigl(ab\,e_3+ux(a+b)+xe_1\bigr),&k\text{ odd}.
 \end{cases}
 \end{gathered}
\end{equation}
The complete five-label coefficients are given in
\eqref{boundary:operatorformula}.
\item All symmetric coefficients satisfy $\overline{f^{(k)}_{R,S}}=f^{(k)}_{R,S}$,
and
\[
 f^{(k)}_{1,-1}=O(u),\qquad
 f^{(k)}_{0,0}=\begin{cases}O(u^2),&k\text{ even},\\O(u),&k\text{ odd}.\end{cases}
\]
The ordinary coefficients satisfy paired simple root-residue cancellation,
and the operator preserves
$\Z[u^{\pm1},x_1^{\pm1},\ldots,x_4^{\pm1}]
[a^{\pm1},b^{\pm1}]^{\mathfrak S_2}$.
\end{enumerate}
\end{theorem}

The boundary correction records the failure of the first-strip
maximal Cohen--Macaulay argument on a deeper mixed-sign support. It is not
determined by the top magnetic coefficient alone.

The extra source factor $x^{-1}$ in Theorem~\ref{thm:main-boundary} is part of
the depth-two top normalization: its electric-zero member is $k=8$, with total
geometric source twist $\ell^8/x^2$ before loop normalization.
The odd-electric correction must be removed as a derived boundary object, not as
a collection of semisimple character summands. This is why the mixed calculation
cannot be obtained by substituting a new electric index in the first-strip formula.
Section~\ref{sec:limits} gives a finite presentation and exact depth test for
general mixed charges, but does not claim that their correction complexes or
full normalized operators have been computed.

The fourth result extends triangular uniqueness beyond the first strip.
Let $\lambda=(r,s)$, $r\geq s$, and let $\mathscr A^{\mathrm{int}}_{\leq\lambda}$
be the operators, normalized by~\eqref{eq:Jglobal} with component units
removed, of integral coherent Coulomb classes represented
in equivariant coherent $G_0$ over $\overline{\mathcal O}_{\lambda}$, with zero
first-node magnetic charge. This is an actual algebra image, not a set defined
only by rational poles and residues. For $D\in\mathscr A^{\mathrm{int}}_{\leq\lambda}$,
write $f_\mu(D)$ for its symmetric coefficient at the dominant label $\mu$.
Put
\[
 \mathscr R=\Z[u^{\pm1},a^{\pm1},b^{\pm1},x_1^{\pm1},\ldots,x_4^{\pm1}].
\]
The bar fixes $a,b,x_v$ and sends $u$ to $u^{-1}$; $\nu_u$ denotes the
Laurent-series valuation at $u=0$ over the fraction field of the other
characters, with $\nu_u(0)=+\infty$.

\begin{theorem}\label{thm:main-uniqueness}
Let $m=4$, $\lambda=(r,s)$ with $r\geq s$, and $p\in\mathscr R$.
Define
\[
 \mathscr L_\lambda(p)=
 \left\{D\in\mathscr A^{\mathrm{int}}_{\leq\lambda}:
 \begin{array}{l}
 f_\lambda(D)=p,\\
 \overline{f_\mu(D)}=f_\mu(D),\quad \nu_u(f_\mu(D))>0
 \\
 \text{for every }\mu<\lambda
 \end{array}\right\}.
\]
Then
\begin{equation}\label{eq:main-unique-cardinality}
 \#\mathscr L_\lambda(p)\leq1.
\end{equation}
Identifying the classes $\mathsf K_k$ of
Theorem~\ref{thm:main-boundary} with their normalized abelianized operators,
one has
\[
 \mathscr L_{(2,-2)}(b^{k-8})=\{\mathsf K_k\}
 \qquad\text{for every }k\in\Z.
\]
\end{theorem}

The order $\mu<\lambda$ is the dominance order in the fixed determinant
component, and the remaining coefficients are determined by Weyl covariance.
The inequality in~\eqref{eq:main-unique-cardinality} is an at-most-one
statement; it does not assert that $\mathscr L_\lambda(p)$ is nonempty for
general $\lambda$ and $p$.

The proof in Theorem~\ref{thm:all-support-uniqueness} uses matter Euler
divisibility on the highest actual incidence stratum to obtain a Laurent
symmetric coefficient. It does not follow from root-residue equations alone:
Remark~\ref{rem:matter-divisibility} gives a pole-free counterexample to that
stronger claim. Thus the theorem supplies neither a converse image theorem
nor a construction of a complete canonical basis.

The paper is organized as follows. Section~\ref{sec:prelim} fixes the coherent,
equivariant, and operator conventions. Section~\ref{sec:geometry} proves the
uniform geometric extension. Section~\ref{sec:localization} derives the shuffle
characters, assembles the full operator, and establishes its structural properties.
Section~\ref{sec:same-sign} treats the same-sign cones and selected first-strip
specializations. Section~\ref{sec:boundary} proves the mixed boundary correction.
Section~\ref{sec:limits} gives the depth criterion for the remaining geometric
problem, proves actual-image triangular uniqueness on arbitrary rank-two
supports, and distinguishes it from residue-only characterizations.

\section{Coherent and operator conventions}\label{sec:prelim}
\subsection{Lattices and the space of triples}
Put $\cO=\C[[z]]$ and $\cK=\C((z))$. The affine Grassmannian
\[\Gr_{\GL_2}=\GL_2(\cK)/\GL_2(\cO)\] parametrizes rank-two $\cO$-lattices in
$\cK^2$. A dominant label $(r,s)$ has $r\geq s$, determinant valuation $r+s$,
and orbit dimension $r-s$. All supports in this paper have their reduced
structures. A magnetic label records a lattice orbit; an electric label is a
representation of its stabilizer and does not change that magnetic orbit.

For the datum~\eqref{eq:datum}, the classical space of triples is
\[
 \cR^{\mathrm{cl}}_{G,\mathbf N_m}
 =\{[g,v]\in G(\cK)\times_{G(\cO)}\mathbf N_m(\cO):gv\in\mathbf N_m(\cO)\}.
\]
We call $v$ the internal matter coordinate and $s=gv$ the external coordinate.
The derived space of triples is the corresponding derived bundle intersection.
On bounded lattice supports it admits the finite presentations used below. The
space of triples is a convolution space, not the Coulomb branch itself.
With the diagonal framing torus $T_F=(\C^\times)^{m-1}$, set
\[
 \widehat G_{\cO}=(G(\cO)\rtimes\C_u^\times)\times T_F,
 \qquad
 \mathcal A_u(G,\mathbf N_m)
 =K^{\widehat G_{\cO}}(\cR^{\mathrm{cl}}_{G,\mathbf N_m}).
\]
The cocharacter $u\mapsto(u^2,u)$ restricts independent loop rotation and matter
dilation: $z$ has weight $u^2$, and matter has dilation weight $u$. Source actions
use inverse source characters. Convolution is composition of lattice modifications;
abelianization identifies it with difference-operator composition~\cite{SS}.

At a finite classical stage $X$ with algebraic acting group $H$, our coherent
convention is
\[
 G_0^H(X)=K_0\bigl(D^b\operatorname{Coh}^H(X)\bigr),\qquad
 [E]=\sum_i(-1)^i[H^i(E)].
\]
The infinite-space notation above means the coherent groups with their prescribed
cutoff transitions. In particular it does not replace coherent K-theory by
$K_0(\operatorname{Perf}X)$ on a singular slice. The reflexive modules used here
need not be locally free. Passing to classical truncation does not change the
K-group, but it changes the coherent category and can discard the module action
needed for a boundary correction. The surrounding coherent pullback formalism is
developed in~\cite{CW25}; our constructions are verified on explicit finite stages.

\begin{remark}\label{rem:involution-comparison}
Here is the explicit comparison with~\cite[Section~3.1]{SSWX}. Set
$F=\C^{m-1}$ and take the disjoint union of arrows
$1\longrightarrow2$ and $\tau2\longrightarrow\tau1$, with the involution
interchanging the two copies. Choose $V_1=\C$, $V_2=V_0$,
$V_{\tau i}=V_i^*$, and only the framing $W_{\tau2}=F^*$ nonzero.
There are no fixed vertices or fixed arrows. Choosing $1\longrightarrow2$
as the positive arrow gives
\[
 \Hom(V_1,V_2)\oplus\Hom(V_{\tau2},W_{\tau2})
 \cong\Hom(\C,V_0)\oplus\Hom(F,V_0)=\mathbf N_m.
\]
The isomorphism is compatible with the gauge and framing actions, with the
framing action on the dual spaces understood contragrediently. Flavor
characters remain equivariant parameters; the framing torus is not gauged.
The parameter denoted $q$ in~\cite{SSWX} equals our $u$: both use loop
weight $q^2=u^2$ and matter dilation $q=u$. This equality precedes the
additional target-character and component normalizations used below; it is
not an assertion that differently normalized displayed operators coincide.

The dominant charge $(0;1,-1)$ is quasi-minuscule, so
\cite[Theorem~2.2]{SSWX} applies to this first member of our first strip.
Its geometric input is a dressed structure sheaf on a resolution and its
proper pushforward. The localization formula determines that class modulo
the coefficient ring, hence modulo a zero-magnetic-charge term. It does
not identify the pushforward with the simple intermediate extension.
For $n>1$, the charge $(0;n,-1)$ is not quasi-minuscule; neither is
$(0;2,-2)$. The saturated dg modules and the derived boundary quotient in
this paper therefore require additional geometric arguments, even though
the ambient Coulomb datum is the same. The homological formula in
\cite[Theorem~2.4]{SSX} is a precursor of this comparison, not a replacement
for its coherent $K$-theoretic or categorical content.
\end{remark}

\subsection{Koszul perversity and localization}
The scalar target cocharacter $\eta$ acts trivially on the magnetic base and with
weight one on each matter column. The base perverse-coherent function is
\[
 p(\cO_\lambda)=-\tfrac12\dim\cO_\lambda.
\]
For an orbit inclusion $i$, the conditions are
$Li^*F\in D^{\leq p(\cO_\lambda)}$ and
$Ri^!F\in D^{\geq p(\cO_\lambda)}$. A pure central weight $v$ is regraded by
$[-v]$. Odd-dimensional components use the half-integral component convention of
Cautis--Williams; no fractional numerical sign is assigned to a shift.
We use the finite-length perverse-coherent heart~\cite{AB} and its Koszul regrading
as in~\cite{CW23}. The finite bundle-intersection lifting theorem applies when the
central cocharacter acts trivially on the base and with weight one on the ambient
bundle. Our lattice supports have the equivariant resolution property from a
linearized Pl\"ucker ample bundle. The t-exact conservative zero-matter detector
will be applied only in this setting, not as a claim about arbitrary derived
pullbacks. 
Coordinate functions have contragredient characters. For an isolated fixed point
on a singular affine slice, choose a smooth ambient representation with function
weights $\chi_i$. The coherent localization of a module is its ambient
K-polynomial divided by $\prod_i(1-\chi_i)$, equivalently its rational Hilbert
character~\cite{Thomason}. It is not the character of its ordinary fiber.
An ambient free resolution suffices; a finite free resolution over the singular
slice is not assumed. The intrinsic Koszul generators are retained in the dg
object even when their Euler factors cancel in a polarization comparison.

\subsection{Characters and normalization}
Let $a,b$ be the target characters. The source specialization is
\[
 (x_1,x_2,\ldots,x_m)=(x_{\mathrm g},z_1,\ldots,z_{m-1}).
\]
The first entry is the dimension-one gauge character and has zero shift throughout.
For $m=4$, the optional reduced framing relation $z_1z_2z_3=1$ makes
$x=x_{\mathrm g}$. All integral identities are proved before specialization.
Use $A,B$ for raw target variables when a simultaneous change to symmetric
variables could otherwise be ambiguous.

For an integral magnetic pair define
\begin{equation}\label{eq:operators}
 \U_{R,S}F(a,b)=F(u^{2R}a,u^{2S}b),\qquad
 B_{R,S}(a,b)=J_{R,S}(a,b)f_{R,S}(u^Ra,u^Sb).
\end{equation}
The operators satisfy
\[
 (B\U_{R,S})(C\U_{P,Q})
 =B(a,b)C(u^{2R}a,u^{2S}b)\U_{R+P,S+Q}.
\]
Thus coefficients do not commute with shifts. The Weyl involution exchanges
$(a,b)$ and $(R,S)$. The word symmetric in~\eqref{eq:operators} denotes the
normalization of $f$, not the assertion that every fixed-label coefficient is
symmetric in $a,b$. Weyl covariance is
$f_{S,R}(b,a)=f_{R,S}(a,b)$. The coefficient bar involution fixes all characters
except $u$, which it sends to $u^{-1}$.
For an integer $d$, the notation $F=O(u^d)$ means that the Laurent expansion
of $F$ at $u=0$, over the fraction field of the remaining independent
characters, has valuation at least $d$. It is not a statement about signs
of coefficients.

For the first-strip family put
\begin{equation}\label{eq:J}
 \begin{aligned}
 V_h(a,b)&=\prod_{l=0}^{h-1}(1-u^{-2l}b/a),\qquad V_0=1,\\
 J_{R,S}(a,b)&=\frac{(-u)^{-S}}{V_{R-S}(a,b)}\quad(R\geq S\geq0),\\
 J_{N-1,-1}(a,b)&=-\frac{u}{V_N(a,b)}
       \prod_{v=1}^m(1+u^{-1}b/x_v).
 \end{aligned}
\end{equation}
The plus signs incorporate simultaneous target negation. All Weyl mates are
obtained by exchanging $a,b$. For the additional four-column sectors in
Sections~\ref{sec:same-sign} and~\ref{sec:boundary}, the calibrated extension is
\begin{equation}\label{eq:Jglobal}
 \begin{aligned}
 P_d(Z)&=\prod_{v=1}^4\prod_{l=1}^d(1+u^{-(2l-1)}Z/x_v),\quad P_0=1,\\
 J_{R,S}(A,B)&=\frac{(-u)^{-S}
 P_{\max(-R,0)}(A)P_{\max(-S,0)}(B)}{V_{R-S}(A,B)}\quad(R\geq S).
 \end{aligned}
\end{equation}
This is the root and matter product of~\cite[(53)--(55)]{Gaiotto}, with the
common determinant-component unit removed. Explicitly, if $\tau^2=-u$, then
$D_{R,S}=\tau^{R+S}J_{R,S}$. The same component unit is inserted only once in an
assembled operator. The factors in~\eqref{eq:Jglobal} are not obtained by
commutatively raising a minuscule operator to a power.

The residue tests considered here use the even ordinary walls $b/a=u^{2h}$.
The partner of $(R,S)$ on this wall is
\begin{equation}\label{eq:generalpartner}
 (R,S)\longmapsto(S+h,R-h).
\end{equation}
For comparison with~\cite[Conjecture 1]{Gaiotto}, we call Weyl covariance and
coefficientwise bar invariance A1, strict lower order $O(u)$ A2, and simple poles
together with the paired even-wall residue identities A3. These elementwise
conditions do not by themselves assert convolution closure or identify an entire
canonical basis.

Distinct walls are considered over the formal generic parameter $u$. Integral
polynomial action means cancellation in the Laurent polynomial ring before any
specialization of $u$; it is not separate evaluation of singular summands at a
root of unity. The extra odd-wall freedom of a larger rational coefficient domain
is addressed in Section~\ref{sec:limits}.

\section{Saturation on the whole magnetic support}\label{sec:geometry}
Every lattice in $S$ satisfies
$z^{N-1}\cO^2\subset L\subset z^{-1}\cO^2$ and has determinant valuation $N-2$.
The finite bundles used in the external polarization are
\begin{equation}\label{eq:bundles}
 V=z^{-1}\cO^2/z^{N-1}\cO^2,\quad
 P=\cO^2/z^{N-1}\cO^2,\quad W=L/z^{N-1}\cO^2,\quad Q=V/W.
\end{equation}
Their ranks are $2N,2(N-1),N,N$. Multiplication by $z$ induces $T:Q\to Q$.
On $U$, set $I=\im(P\to Q)=TQ$; it has rank $N-1$. The electric line $\ell$
is the leading-coefficient image of $L$ in the constant representation $V_0$.

\subsection{The reduced framed model}
\begin{lemma}\label{lem:singularity}
The variety $S$ is normal, Gorenstein, and has rational singularities. In particular
it has canonical, hence klt, singularities. Its boundary has codimension two.
\end{lemma}
\begin{proof}
Let $\cN_N\subset\mathfrak{gl}_N$ be the reduced nilpotent cone. Framing $Q$ gives
an isomorphism
\begin{equation}\label{eq:frames}
 \Fr(Q)\simeq
 \{(T,v_1,v_2)\in\cN_N\times(\C^N)^2:
              \C[T]v_1+\C[T]v_2=\C^N\}.
\end{equation}
A frame sends the images of $z^{-1}e_1,z^{-1}e_2$ to $v_1,v_2$. Conversely, the
generating pair defines a surjection $(\C[z]/z^N)^2\to\C^N$. Its kernel is a
rank-$N$ subbundle, stable under $z$, and reconstructs $W$ and $L$. These constructions
are algebraic and inverse. Cayley--Hamilton supplies $T^N=0$ on $\cN_N$.

The right side is an open subset of the product in~\eqref{eq:frames}, since generation
is a matrix-rank condition. Its dimension is $N^2+N$, and it is irreducible. The quotient
by the frame group has dimension $N$, lies in the stated lattice bounds and determinant
component, and contains $U$. Thus it is precisely the reduced closure $S$. The reduced
qualification matters: the equations $T^N=0$ alone need not impose the nilpotent-cone
scheme over nonreduced test rings.

The nilpotent cone is normal and a complete intersection. It has rational singularities
by the vanishing for its Springer resolution; see Hesselink~\cite[Theorem A]{Hesselink}. These properties
pass to the open product and descend along the Zariski-locally trivial principal
$\GL_N$-bundle $\Fr(Q)\to S$. A rational Gorenstein complex variety is canonical.
Finally, dominance gives exactly the strata in~\eqref{eq:support}; their dimensions
$N-2j$ prove the codimension assertion.
\end{proof}

\subsection{The electric line and its actual cyclic cover}\label{sec:electric}
On $U$, the cokernel $C=Q/TQ$ is a line. It is the cokernel, not the kernel of $T$,
that enters the leading-coefficient sequence:
\begin{equation}\label{eq:leading}
 0\longrightarrow t^{-1}\ell\longrightarrow t^{-1}V_0
       \longrightarrow C\longrightarrow0.
\end{equation}
The regular-nilpotent filtration has successive quotients $Ct^a$, $0\leq a<N$.
Consequently
\begin{equation}\label{eq:electricpower}
 \ell=t^{-1}(\det V_0)C^{-1},\qquad
 \ell^N=E|_U,\qquad
 E=(\det V_0)^N(\det Q)^{-1}t^{N(N-3)/2}.
\end{equation}
The bundle $E$ exists on all of $S$; no global root or trivialization is assumed.
At the fixed point $\lambda_j$, the two quotient chains have weights
$At^{-1},\ldots,At^{N-j-2}$ and $Bt^{-1},\ldots,Bt^{j-2}$. Hence
\begin{equation}\label{eq:detQ}
 \det Q=A^{N-j}B^jt^{N(N-3)/2-j(N-j)},\qquad
 E_j=A^jB^{N-j}t^{j(N-j)}.
\end{equation}

\begin{lemma}\label{lem:cyclicmcm}
Let $X$ be normal and klt over $\C$, let $U_X\subset X$ have complement of codimension
at least two, and let $L$ be a line on $U_X$ with $L^N=E_X|_{U_X}$ for a line $E_X$
on $X$. Then the ordinary extensions $\iota_*L^k$ are coherent and maximal
Cohen--Macaulay for every $k\in\Z$.
\end{lemma}
\begin{proof}
Locally trivialize $E_X$ and form the finite étale cyclic torsor on $U_X$ with
algebra $\bigoplus_{r=0}^{N-1}L^r$. Normalize $X$ in its finite generic algebra,
component by component if the torsor is disconnected. Excellence gives a finite
normal scheme $f:X'\to X$. Its restriction to $U_X$ is the given torsor, and the
inverse image of the boundary still has codimension at least two. Thus $f$ is
quasi-étale. Finite quasi-étale covers preserve klt singularities
\cite[Proposition 5.20]{KM}; see also~\cite[Reminder 2.14]{GGK}. Therefore $X'$ is
Cohen--Macaulay.

A parameter system in $\cO_{X,x}$ is a parameter system in every local ring above it:
finiteness and the dimension formula give the same local dimension. It is regular
there, so $(f_*\cO_{X'})_x$ is an MCM module over $\cO_{X,x}$. Deck projectors split
this finite module into coherent MCM summands. Each is torsion-free and $S_2$, hence
reflexive, and is the unique reflexive extension of its line on $U_X$. This identifies
the summands with $\iota_*L^r$. Writing $k=qN+r$ and tensoring with $E_X^q$ proves
the assertion for every integer $k$. Flatness of the finite cover is not needed.
\end{proof}

The local cover also fixes the electric sign. For a cyclic vector $c$ define
\[
 D(c)=\det(c,Tc,\ldots,T^{N-1}c)
\]
and adjoin $y$ with $y^N=D(c)$. We use the deck
action $y\mapsto\zeta y$. Replacing $c$ by $f(T)c$ changes $y$ by $f(0)y$, so
$c/y$ is a normalized frame of $C$ of deck character $-1$. The reciprocal electric
frame $t^{-1}AB/(c/y)$ has deck character $+1$. If $c_j^N=E_j$, invariant descent gives
\begin{equation}\label{eq:deckelectric}
 \ch(\ell^{[k]})=c_j^kH_{-k}^{N,j}(t,B/A),\qquad k\in\Z.
\end{equation}
Here brackets denote reflexive extension, not a cohomological shift. Equation
\eqref{eq:deckelectric} retains the full power $c_j^k$: increasing $k$ by $N$ tensors
by $E$, rather than returning to the same line.

\subsection{Matter terms and the full dg module}
Over $Y=S\times P^m$ the finite derived incidence and its algebra are
\begin{equation}\label{eq:dg}
 \cD=W^{\oplus m}\mathbin{\times^R}_{S\times V^m}Y,\qquad
 \cB=(\cO_Y\otimes\Lambda((Q^*)^{\oplus m}),\iota_\sigma).
\end{equation}
There are $mN$ ambient generators, all in degree $-1$. In this exterior notation the
$m$ dual copies carry the respective factors $u^{-1}x_r$; the same convention is
used for $I^*$ below. On $Y_U$, the intrinsic complex is
\[
 K_k=(\cO_{Y_U}\otimes\Lambda((I^*)^{\oplus m})\otimes\ell^k/x,\iota_\sigma),
\]
with $m(N-1)$ intrinsic generators; $Q^*\to I^*$ gives its $\cB$-action. Put
\begin{equation}\label{eq:terms}
 T_l=\iota_*\bigl(\Lambda^l((I^*)^{\oplus m})\otimes\ell^k/x\bigr),
 \qquad 0\leq l\leq m(N-1).
\end{equation}

\begin{proposition}\label{prop:extension}
Every $T_l$ is coherent MCM. Ordinary termwise pushforward of $K_k$ gives a bounded
coherent dg $\cB$-module $\cM_k$ with $\cM_k^{-l}=\pr_S^*T_l$. Its differential,
ambient action, and positive-loop linearizations extend uniquely.
\end{proposition}
\begin{proof}
At every boundary point the Jordan type of $T$ is $(N-j,j)$, so its rank is exactly
$N-2$, independently of $j>0$. On a neighborhood choose an invertible $(N-2)$-minor.
Gaussian elimination gives $Q=F\oplus H$, of ranks $N-2$ and $2$, such that on $U$
one has $I=F\oplus J$, with $J$ a line in $H$. Since $H/J=C$,
\[
 J^*=C(\det H)^{-1}
     =\ell^{-1}t^{-1}(\det V_0)(\det H)^{-1}.
\]
Exterior expansion expresses every term in~\eqref{eq:terms} locally as a finite sum
of $\ell^{k-r}$ tensored with an extending locally free sheaf, $0\leq r\leq m$.
All $m(N-2)$ free exterior directions remain in this expansion. Lemma~\ref{lem:cyclicmcm}
proves coherence and MCM at the boundary, and on $U$ the assertion is immediate.

Polynomial pullback of $T_l$ is reflexive and agrees with the open term. Reflexive
Hartogs extension across $Y\setminus Y_U$ identifies it with ordinary pushforward
\cite[Tags 0AVB and 0AVC]{Stacks}.
For any sheaf $A$ and open sheaf $B$, the adjunction
$\Hom_Y(A,\iota_*B)=\Hom_{Y_U}(\iota^*A,B)$ extends every differential and wedge-action
map. Taking $A$ to be a term, or its tensor product with a locally free term of $\cB$,
gives the required maps. The same adjunction proves the differential-square, Leibniz,
associativity, and unit identities. No tensor product of reflexive sheaves is assumed
to be reflexive. A nontrivial electric twist gives a module, not an intrinsic algebra.

At any finite cutoff let $H_0$ be a sufficiently large smooth target/source jet group,
including loop and flavor tori. On $H_0\times Y$, the action and projection pullbacks
of the terms are reflexive because those maps are smooth. The intrinsic open
linearization and its inverse extend uniquely. On $H_0\times H_0\times Y$, uniqueness
proves the cocycle, and likewise compatibility with the differential and $\cB$.
Thus the action is algebraic, not merely infinitesimal. Local cyclic covers and local
Gaussian elimination were needed only for the underlying MCM property; no equivariant
choice of either is required.
\end{proof}

\subsection{Depth and the intermediate extension}
\begin{proposition}\label{prop:simple}
In the Cautis--Williams Koszul-perverse category, $\cM_k[N/2-k]$ is the simple
intermediate extension of its open incidence line object. It defines a bounded-support
simple in the full finite-stage category.
\end{proposition}
\begin{proof}
Set $d=N/2$. Remove strata smaller than $S_j$ so that $S_j$ is closed, with inclusion
$i_j$ and ideal of height $2j$. For each MCM term, the grade of this ideal is $2j$.
The local depth formula for grade~\cite[Tag 0AVY]{Stacks} gives
\begin{equation}\label{eq:depth}
 Li_j^*T_l\in D^{\leq0},\qquad Ri_j^!T_l\in D^{\geq2j}.
\end{equation}
The first bound holds for any sheaf; it does not assert bounded negative Tor. The
second follows from the vanishing of $\Ext^q(\cO_{S_j},T_l)$ for $q<2j$. For middle
perversity $p(S_j)=j-d$, the shifted bounds are strictly on either side of $p(S_j)$
when $j>0$. Thus $T_l[d]$ is perverse coherent and has no boundary-supported subobject
or quotient. Indeed, take the largest stratum in such a support and use adjunction
with the strict bounds after removing smaller strata; induction over the finite
stratification completes the argument, also for nonreduced supports.

The target scalar cocharacter $\eta$ acts trivially on $S$, with weight one on each
matter column. The weight of $T_l$ is $k-l$. The linearized Pl\"ucker ample line on
$S\subset\operatorname{Gr}(N,V)$ gives the equivariant resolution property; its
restrictions supply the same property on the invariant opens used above. The standard
finite-length perverse-coherent heart and its $\eta$-regrading therefore satisfy the
hypotheses of the Cautis--Williams Koszul t-structure lifting and
weight-regrading constructions~\cite{CW23}.
The lifted detector $\Phi=L\sigma_0^*i_*$, for the matter zero section and the derived
closed embedding, is t-exact and conservative in this setting. The matter-coordinate
sequence is regular on each polynomial term; hence
\begin{equation}\label{eq:detector}
 \Phi(\cM_k)=\bigoplus_{l=0}^{m(N-1)}T_l[l].
\end{equation}
The differential becomes zero at the zero section. Regrading a weight-$v$ object by
$[-v]$ sends $T_l[d]$ to $T_l[l+d-k]$. Thus~\eqref{eq:detector} after $[d-k]$ lies
in the base Koszul heart with no boundary subquotients. T-exact conservativity proves
the same two assertions for $\cM_k[d-k]$; support compatibility excludes boundary
subobjects and quotients.

On $U$ the independent incidence equations resolve the classical incidence vector
bundle. The open line is a one-dimensional irreducible stabilizer representation,
so the Cautis--Williams simple-object correspondence~\cite{CW23} identifies the
shifted open object as simple. A proper
subobject of its extension either restricts to zero or has quotient restricting to
zero. The boundary exclusions prove simplicity and the intermediate-extension property.

For completeness the cutoff transitions preserve this object, rather than adding
new exterior factors. A positive matter cutoff is a polynomial pullback: equations
depend only on $P$ because $z^{N-1}\cO^2\subset L$. A large common jet order is
$\max(N,M+1)$ at cutoff $P_M=\cO^2/z^{M+1}\cO^2$; the compatible tower need not
factor through one fixed jet group. Enlarging the negative ambient bound gives
$Q_R^*\twoheadrightarrow Q^*$ and a dg map $\cB_R\to\cB$. The transition is
restriction of scalars, equivalently derived closed pushforward, not tensoring with
an extra exterior algebra. Enlarging the magnetic support is closed pushforward.
These are the simple-preserving transitions of the Cautis--Williams finite-stage
construction~\cite{CW23}.

For odd $N$, the symbol $[d-k]$ belongs to the half-integrally graded determinant
component of the Cautis--Williams grading convention~\cite{CW23}. It is not a fractional shift in an ordinary
unmodified derived category, and will not be used to assign a fractional numerical sign.
\end{proof}

\section{Uniform characters and operator localization}\label{sec:localization}
\subsection{The characters of two row normalizations}\label{sec:characters}
Fix $1\leq q=j\leq\lfloor N/2\rfloor$, put $p=N-q$, and take the block Slodowy
slice $\Sigma=e+\ker(\operatorname{ad}f)\subset\mathfrak{gl}_N$ for Jordan blocks
$(p,q)$. Write $X=\Sigma\cap\cN_N$. Its dimension is $2q$. We use positive
Kazhdan \emph{function grading}: equivalently, use the inverse point action when
characters of functions are defined by inverse pullback. This convention gives the
weights below and is maintained throughout.

The framed model~\eqref{eq:frames} gives the transverse comparison with the lattice
support: taking a local section of the frame bundle and then a transverse slice to
the nilpotent conjugacy orbit identifies the completed local germ of $S$ at
$\lambda_j$ with its smooth orbit directions times $X$. The quotient bundle and
its cyclic-vector determinant pull back under this comparison. Thus the local
electric cover is precisely $\widetilde X$, with the chain characters below;
its Hilbert localization is multiplied by the smooth orbit Euler factor in
Section~\ref{sec:assembly}.

The diagonal function weights on $\Sigma$ are $t,\ldots,t^p$ and $t,\ldots,t^q$.
With two quotient-chain characters
\begin{equation}\label{eq:chains}
 A_r=at^{r-1}\ (0\leq r<p),\qquad B_s=bt^{s-1}\ (0\leq s<q),\qquad \rho=b/a,
\end{equation}
the off-diagonal weights are $t^l\rho$ and $t^{p-q+l}/\rho$, $1\leq l\leq q$.
Indeed the $\mathfrak{sl}_2$ decomposition of the two off-diagonal Hom spaces has
one lowest-weight direction in each of $q$ summands. In balanced coordinates
$\rho=t^{(p-q)/2}R$, both lists become $t^{(p-q)/2+l}R^{\pm1}$, which have positive
scalar degree. Thus the weights describe the actual Slodowy slice, not an assumed
identification with an arbitrary companion matrix slice.

Let $\widetilde X$ be the normalization in the pullback of the universal regular-orbit
cover. Its deck summands are $M_i$, indexed by $y\mapsto\zeta y$, with Hilbert
characters $H_i$. Define
\begin{equation}\label{eq:shuffledata}
 \begin{gathered}
 Q_{\det}=a^pb^qt^{d_{p,q}},\qquad
 d_{p,q}=\frac{p(p-3)+q(q-3)}2,\\
 D_{p,q}(t,\rho)=\prod_{l=1}^q(1-t^l\rho)(1-t^{p-q+l}/\rho),\\
 I_i=\{k:1\leq k\leq N,\ \lceil ik/N\rceil>\lceil i(k-1)/N\rceil\},\\
 h_0=0,\qquad h_i=\frac{N-\gcd(i,N)}2\quad(0<i<N).
 \end{gathered}
\end{equation}
Let $\Sh(p,q)$ consist of the interleavings of two descending lists,
\[
 (A_{p-1},\ldots,A_0)\quad\text{and}\quad(B_{q-1},\ldots,B_0).
\]
For a word $z=(z_1,\ldots,z_N)$,
put
\[
 C(z)=\prod_{\substack{k<l\\z_k,z_l\text{ in different lists}}}
                  \frac{1-tz_k/z_l}{1-z_k/z_l}.
\]

\begin{theorem}\label{thm:shuffle}
The character of the actual saturated summand $M_i$ is
\begin{equation}\label{eq:shuffle}
 H_i^{N,j}(t,\rho)=
 \frac{t^{-h_i}Q_{\det}^{-i/N}}{D_{p,q}(t,\rho)}
       \sum_{z\in\Sh(p,q)}C(z)\prod_{k\in I_i}z_k,
 \qquad 0\leq i<N.
\end{equation}
There are exactly $\binom Nq$ summands. The result is independent of the common
scale of $a,b$. Its numerator $P_i=D_{p,q}H_i$ is a finite Laurent character with
nonnegative integer multiplicities on the genuine finite-cover character lattice.
Its total multiplicity is $\binom Nq$, and
\begin{equation}\label{eq:H0}
 H_0(t,\rho)=\frac{\genfrac{[}{]}{0pt}{}{N}{q}_t}{D_{p,q}(t,\rho)},\qquad
 \genfrac{[}{]}{0pt}{}{N}{q}_t=\frac{\prod_{l=1}^N(1-t^l)}
 {\prod_{l=1}^p(1-t^l)\prod_{l=1}^q(1-t^l)}.
\end{equation}
\end{theorem}

\subsubsection{The graded Graham input and the composite sectors}
We first explain why the Weyl expression computes a coherent normalization summand,
not a virtual fit. Let $\alpha_r=\eps_r-\eps_{r+1}$ and set
\begin{equation}\label{eq:lambdai}
 a_r=\frac{(-ir)\bmod N}{N},\qquad
 \lambda_i=\sum_{r=1}^{N-1}a_r\alpha_r.
\end{equation}
Graham's canonical-sheaf pushforward, vanishing, and normalization results
\cite{Graham}
apply to all weight/root-lattice cosets, including those on walls. They identify the
normalization summands through the Springer map $\mu:T^*\mathcal F\to\cN_N$ and
the projection $\pi:T^*\mathcal F\to\mathcal F$.

The scalar refinement is as follows. On a simple-root toric coordinate, the base
canonical monomial has exponent $1$; in sector $i$ it has exponent $1-a_r$.
In particular a zero residue contributes exponent $1$, not $0$. Their degree ratio
is $t^{-\sum a_r}$. All other root directions and flag canonical factors are unchanged.
Dividing the canonical identities by the same pulled-back symplectic volume gives
the graded identity
\begin{equation}\label{eq:gradedGraham}
 M_i^{\mathrm{global}}\simeq
 t^{-h_i}R\mu_*\pi^*\cL_{\lambda_i},
\end{equation}
with the right side concentrated in degree zero. A finite scalar-torus lift permits
the displayed fractional weights; it is not an arbitrary choice of grading shift.

For $a_0=a_N=0$, one has $a_r=\lceil ir/N\rceil-ir/N$. Therefore
\begin{equation}\label{eq:ceil}
 \lambda_i=\sum_{k\in I_i}\eps_k-\frac{i}{N}\sum_{k=1}^N\eps_k,
 \qquad |I_i|=i.
\end{equation}
If $g=\gcd(i,N)$, the nonzero residues among $-ir\bmod N$ are $g$ copies of
$g,2g,\ldots,N-g$. Their sum divided by $N$ is $(N-g)/2$. This proves
$\sum a_r=h_i$, including all nonprimitive sectors of composite $N$.

Flag localization of~\eqref{eq:gradedGraham}, with the ambient coordinate denominator
already multiplied in, yields the $\mathfrak{gl}_N$ K-polynomial
\begin{equation}\label{eq:Weyl}
 \begin{aligned}
 K_i={}&t^{-h_i}(1-t)^N\left(\prod_{r=1}^Nv_r\right)^{-i/N}\\
 &\quad\times\sum_{w\in\mathfrak S_N}
 \left(\prod_{k\in I_i}v_{w(k)}\right)
 \prod_{k<l}\frac{1-tv_{w(k)}/v_{w(l)}}{1-v_{w(k)}/v_{w(l)}}.
 \end{aligned}
\end{equation}
The auxiliary $v_r$ in this equation are ambient eigenvalue characters, not the $m$
source characters. The factor $(1-t)^N$ includes the scalar direction of
$\mathfrak{gl}_N$; it is important when passing from the semisimple Springer setting
to the indicated ambient K-polynomial.

\subsubsection{Restriction and the exact shuffle reduction}
We justify ordinary restriction before evaluating~\eqref{eq:Weyl}. The action map
$\GL_N\times\Sigma\to\mathfrak{gl}_N$ is smooth along $X$. The pullback of the
global normalization is normal under this smooth map; equivariance identifies it
with $\GL_N$ times the ordinary slice pullback. Normality descends. Finiteness and
agreement with the torsor over the regular locus therefore identify the ordinary
pullback with $\widetilde X$, componentwise if necessary.

Graham's normalization theorem~\cite{Graham} makes the global normalization
Cohen--Macaulay, and its
deck summands are MCM on $\cN_N$. The affine slice has dimension $N+2q$, whereas
$\dim\cN_N=N^2-N$ and $\dim X=2q$. Its $N^2-N-2q$ affine-linear equations thus
have exactly the expected codimension on every nonzero summand. On an MCM module
they form a regular sequence. Hence the restricted ambient resolution has no higher
Tor. Dividing its K-polynomial by the smooth slice Euler denominator really gives
$\ch M_i$. Equivariant localization here is the singular coherent localization of
Thomason~\cite{Thomason}, not a fiber evaluation of a nonfree module.

Now substitute the two chains~\eqref{eq:chains} in~\eqref{eq:Weyl}. Same-chain
denominators are nontrivial powers of $t$ and cross-chain denominators retain the
independent ratio $b/a$, so none is identically zero in $\C(t,a,b)$. If $A_r$ precedes
$A_{r+1}$, the numerator contains $1-tA_r/A_{r+1}=0$. Thus every surviving permutation
orders each chain descending. Conversely, no descending shuffle has an identically
zero numerator. Precisely $\binom Nq$ terms remain; no equal-weight substitution
inside a singular summand is used.

For a descending chain of length $d$, the internal factors telescope:
\begin{equation}\label{eq:telescope}
 \prod_{h=1}^{d-1}\left(\frac{1-t^{h+1}}{1-t^h}\right)^{d-h}
       =\frac{\prod_{l=1}^d(1-t^l)}{(1-t)^d}.
\end{equation}
After multiplying by $(1-t)^N$, the two internal contributions are
$(t;t)_p(t;t)_q$, where
\[ (t;t)_d=\prod_{l=1}^d(1-t^l). \]
The affine slice denominator is exactly $(t;t)_p(t;t)_qD_{p,q}$.
Cancellation proves~\eqref{eq:shuffle}. Equation~\eqref{eq:ceil} proves independence
of the common chain scale and also shows why the full scalar part of
$Q_{\det}^{-i/N}$ must be retained.

\subsubsection{A parameter ring and positive finite numerators}
Retain all $2q$ off-diagonal slice coordinates. Setting them to zero leaves the two
regular nilpotent Kostant sections. Their characteristic coefficients form polynomial
coordinates: the coefficient of weight $k$ has nonzero linear term in the coordinate
of $f_m^{k-1}$, since $\operatorname{tr}(e_m^{k-1}f_m^{k-1})\ne0$, plus a polynomial
in smaller-weight coordinates. Thus the change to characteristic coefficients is
weighted triangular.

Write the resulting monic block polynomials as $A(z),B(z)$, of degrees $p,q$.
Nilpotency says $A(z)B(z)=z^N$, so both factors have all roots zero. The reduced
common zero locus of the off-diagonal coordinates in $X$ is therefore the origin.
The positive centered grading shows that these coordinates form a homogeneous
system of parameters. Their polynomial ring $A_0$ has dimension $2q$, and
$\cO(X)$, the normalization, and every $M_i$ are finite over $A_0$.

The MCM property makes this parameter sequence regular on $M_i$. A finite graded
MCM module over a polynomial ring is graded free: lift a homogeneous basis of its
Artinian quotient and use regularity, or apply the graded Auslander--Buchsbaum
argument. Hence
\begin{equation}\label{eq:free}
 M_i\simeq\bigoplus_\nu A_0m_{i,\nu},\qquad
 P_i=D_{p,q}H_i=\sum_\nu\wt(m_{i,\nu}).
\end{equation}
This proves cancellation of apparent cross-ratio poles and nonnegative integral
multiplicities without interpreting the signs of individual shuffle summands.
For $i=0$ the Artinian quotient is the complete intersection with variable weights
$1,\ldots,p$ and $1,\ldots,q$ and relation weights $1,\ldots,N$. Its Hilbert series
is $\genfrac{[}{]}{0pt}{}{N}{q}_t$ and its length is $\binom Nq$. Each sector is a line on the dense
regular locus, hence has rank one over the base; its $A_0$-rank is also $\binom Nq$.
This proves all remaining assertions of Theorem~\ref{thm:shuffle}.

\subsection{Matter polarization and the normalized coefficients}\label{sec:assembly}
We now identify~\eqref{eq:mainformula} with the actual class, before imposing any
root-residue equation. Use raw target characters $A,B$ in this section, to distinguish
the ordinary coordinates from the symmetric coordinates $a,b$.

At a lower fixed point $\lambda_j$, both magnetic entries are nonnegative. The two
Jordan chains provide a constant invertible $(N-2)$-minor. In the positive balanced
grading that minor has degree zero and value one, hence is identically one on the
slice. Equivariant elimination supplies the splitting $I=F\oplus J$ there, where
$F$ consists of the $N-2$ nonnegative quotient modes and the residual target $H$ has
weights $A/t,B/t$. From~\eqref{eq:leading} and $\det H=ABt^{-2}$ one obtains
\begin{equation}\label{eq:residualJ}
 J^*=C(\det H)^{-1}=t\ell^{-1}.
\end{equation}
Each source-function multiplier is $u^{-1}x_r$. Thus a residual exterior generator
has compensation $ux_r/c_j$, where $c_j^N=E_j$, and the $r$-th residual exterior
power coupled to $\ell^k$ uses sector $r-k$.

The free generators have weights $u^{-1}x_r/F$. They are precisely the external
constrained modes $0,\ldots,\lambda_1-1$ and $0,\ldots,\lambda_2-1$.
Because $W\subset P$ at a lower point, there are no internal constrained modes.
In a common smooth ambient vector bundle, comparison of the external and internal
zero sections cancels the Euler factors of these $m(N-2)$ free directions. This
comparison is an equality for every coherent K-class, computable through an ambient
free resolution; it requires neither deleting generators from $\cM_k$ nor a finite
free resolution over the singular slice. The smooth orbit contributes
$1/V_{N-2j}(A,B)$.

Simultaneous target negation changes the sign of every matter generator, cancelling
the alternating sign of its exterior Euler class. The sum of residual contributions
therefore uses $+e_r$. The character $H_i$ has target-central weight zero, while
$c_j$ has weight one. This observation fixes the same negation on the electric frame.

It remains to specify the single scalar. The determinant component is $N-2$, with
baseline $[(N-2)/2]$. The simple has shift $[N/2-k]$, whose difference from the
baseline is the integer $1-k$. Its relative parity is $(-1)^{1-k}$. Target negation
on $\ell^k$ contributes $(-1)^k$, giving the constant relative sign minus. This
uses no fractional numerical sign when $N$ is odd.

At the top point the bare internal coefficient after negation is
\begin{equation}\label{eq:topmatch}
 K_{\mathrm{int}}=\frac{\prod_r(1+ux_r/B)}{V_N(A,B)},\qquad
 \frac{J_{N-1,-1}}{K_{\mathrm{int}}}=-u^{1-m}B^m/x.
\end{equation}
The desired leading symmetric function $b^{k-m}$ becomes
$(u^{-1}B)^{k-m}$ in ordinary coordinates. Consequently the already source-normalized
object $\ell^k/x$ is multiplied, in the rephased ordinary convention, by
$-u^{1-k}$. If one begins with $\ell^k$, the total multiplier is
$-u^{1-k}/x$. In particular $x^{-1}$ appears exactly once.

Combining these facts with~\eqref{eq:deckelectric} gives the actual lower ordinary
localization
\begin{equation}\label{eq:rawlower}
 B_{\lambda_j}^{(k)}(A,B)=
 -\frac{u^{1-k}c_j^k}{xV_{N-2j}(A,B)}
       \sum_{r=0}^me_r(u/c_j)^rH_{r-k}^{N,j}(u^2,B/A).
\end{equation}
This calculation uses the saturated terms already constructed. In particular the
Hilbert character is not the character of a naive root hypersurface.

Divide~\eqref{eq:rawlower} by
$J_{\lambda_j}=(-u)^{-(j-1)}/V_{N-2j}$ and substitute
\begin{equation}\label{eq:symshift}
 A=u^{-(N-1-j)}a,\qquad B=u^{-(j-1)}b.
\end{equation}
Equation~\eqref{eq:detQ} gives, exactly,
$c_j\mapsto u\delta_j$ and $B/A\mapsto u^{N-2j}b/a$.
The remaining sign is $(-1)^j$, and the scalar power is $u^j$.
This proves~\eqref{eq:mainformula}. The top term was separately computed in
\eqref{eq:topmatch}; extending the lower formula to $j=0$ would miss its changed
matter polarization. The common factor $\tau^{N-2}$, if desired, is inserted only
once after the ordinary operator has been assembled.

\subsection{Balanced duality and strict lower order}\label{sec:bar}
For $d_j=p-q=N-2j$, set
\begin{equation}\label{eq:G}
 G_i(u,\rho)=H_i^{N,j}(u^2,u^{d_j}\rho).
\end{equation}
The off-diagonal function weights are now
$u^{d_j+2l}\rho^{\pm1}$, $1\leq l\leq j$, and the diagonal weights are positive even
powers of $u$. The lower coefficient in~\eqref{eq:mainformula} becomes
\[
 \frac{(-1)^ju^j}{x}\sum_{r=0}^me_r\delta_j^{k-r}G_{r-k}(u,b/a).
\]

\begin{lemma}\label{lem:duality}
The actual cover-sector characters satisfy
\begin{equation}\label{eq:threeidentities}
 \begin{aligned}
 G_i(u,\rho)&=G_{-i}(u,\rho^{-1}),\\
 G_i(u^{-1},\rho^{-1})&=u^{2j}G_{-i}(u,\rho),\\
 G_i(u^{-1},\rho)&=u^{2j}G_i(u,\rho).
 \end{aligned}
\end{equation}
\end{lemma}
\begin{proof}
Write $V=V_p\oplus V_q$ as a sum of irreducible $\mathfrak{sl}_2$ modules. Choose
invariant nondegenerate forms with block matrix $J$, using the inverse form matrix
if necessary. The contragredient automorphism
$\theta(X)=-JX^{\mathsf T}J^{-1}$ fixes $e,h,f$, hence preserves the Slodowy slice,
its nilpotent intersection, and its regular locus. It commutes with the balanced
Kazhdan action (and its inverse), so preserves the positive function grading. On
block scalars the group automorphism $g\mapsto Jg^{-\mathsf T}J^{-1}$ inverts both
characters; therefore it inverts the balanced flavor ratio $\rho$.

The regular universal cover is $\SL_N/C^0$, with $C$ the regular-nilpotent centralizer.
After conjugating the reference regular nilpotent if needed, $\theta$ acts on this
cover. Since $\theta(\zeta I)=\zeta^{-1}I$, it reverses deck characters. Pull back
to the slice and extend uniquely to the finite normalization. Compatibility with
the lifted connected tori follows from the homogeneous-space construction, or from
uniqueness: two normalized lifts differ by a map from a connected torus to the finite
deck group. This proves the first identity. The forms on $V_p,V_q$ may have different
symmetry types; then $\theta^2$ is an inner block-scalar automorphism. The argument
does not require $\theta$ to square to the identity.

The slice has dimension $2j$ and a symplectic volume of character $u^{2j}$ in the
function/differential grading, with trivial flavor character. Its regular-nilpotent
complement has codimension at least two. The normalized cover is finite quasi-étale
and CM; the pulled-back volume generates its reflexive canonical sheaf and has
trivial deck character. Thus it is Gorenstein with this canonical generator.
Multigraded Cohen--Macaulay reciprocity, retaining the finite deck grading, gives the
second identity: inversion sends sector $i$ to $-i$, the canonical factor is $u^{2j}$,
and the dimension sign is positive. Substitute $\rho^{-1}$ for $\rho$ in that identity
and use the first to obtain the third. Ordinary simultaneous-variable reciprocity
alone would not imply this fixed-sector bar identity.
\end{proof}

\begin{lemma}\label{lem:degree}
Each $G_i$ has nonnegative integral $u$-degrees, and every weight in sector $i$ has
flavor exponent congruent to $-ji/N$ modulo $\Z$. Hence
$\delta_j^{k-r}G_{r-k}(u,b/a)$ has only integral target characters for every $k\in\Z$.
\end{lemma}
\begin{proof}
The coordinate ring of $X$ is positively graded with degree-zero part $\C$. The finite
normalization is integral over it. Even allowing an initial finite grading cover,
a nonzero homogeneous element of negative degree cannot satisfy a monic integral
equation: its leading power is the unique term of smallest degree, since every
coefficient has nonnegative degree. It would be nilpotent, contrary to reducedness.
Thus all normalized degrees are nonnegative.

To check that scalar degrees are integers, take the first quotient-chain vector
$c=z^{-1}e_1$ on its cyclic open chart. Its determinant has weight
\begin{equation}\label{eq:rootweight}
 \wt D(c)=r^{-j}t^{j(N-j)},\qquad
 r=u^{N-2j}\rho,\qquad
 \wt D(c)=u^{jN}\rho^{-j}.
\end{equation}
Therefore the root $y$ with $y^N=D(c)$ has weight $u^j\rho^{-j/N}$. The generic torsor algebra and its
normalization carry an honest integral scalar grading; only the flavor torus needs
a finite cover. Since sector $i$ is generically $\C(X)y^i$, its flavor exponents
have residue $-ji/N$. Write $\delta_j=b\rho^{-j/N}$. For $i\equiv r-k\pmod N$,
the fractional exponent of $\delta_j^{k-r}G_i$ is
$-j(k-r+i)/N\in\Z$. Its remaining $b$-power is also integral. This argument works
for negative $k$ and nonprimitive sectors alike.
\end{proof}

By Lemma~\ref{lem:duality}, $u^jG_i$ is fixed by $u\mapsto u^{-1}$ with flavor fixed.
The other factors of each summand in~\eqref{eq:mainformula} are $u$-independent.
Lemma~\ref{lem:degree} gives $u^jG_i=O(u^j)$; the denominator has constant term one
at $u=0$. This proves bar invariance and strict lower order term by term, not by
assuming positivity of a fitted rational numerator. Weyl covariance comes from
equivariance of the geometric class. When $p=q$, the balanced flavor inversion and
the compensating character also make the central coefficient symmetric, so it is
counted once.

A useful consequence for the finite numerators is
\begin{equation}\label{eq:preciprocity}
 H_i(t^{-1},\rho^{-1})=t^jH_{-i}(t,\rho),\qquad
 P_i(t^{-1},\rho^{-1})=t^{-pq}P_{-i}(t,\rho).
\end{equation}
The second identity uses that the product of the $2q$ denominator weights is
$t^{q(p+1)}$. In centered coordinates the numerator degrees lie between $0$ and
$pq$ in $t$-degree. The lower bound uses integrality over the positive base, and
the upper bound then follows from reciprocity; reciprocity alone would not establish
the lower bound.

\subsection{Root residues and integral polynomial action}\label{sec:residues}
The final step is a geometric image theorem with its conventions matched explicitly.
It would be invalid to infer the required residue law from an arbitrary collection
of rational local characters merely because they look geometric. Sections
\ref{sec:geometry}--\ref{sec:assembly} supplied both a genuine coherent class and
its actual abelianized coefficients.

\subsubsection{The image theorem and the wall partner}
In Schrader--Shapiro~\cite[Section 2.3, equation (2.6)]{SS}, integral abelianization
embeds the quiver Coulomb algebra in the difference algebra localized at even-root
binomials. Definition 2.5 imposes simple poles and paired root residues; Proposition
2.7(1),(3) identifies the corresponding pure-$\GL_d$ residue algebra. Crucially,
equation (2.20), at the start of Section 2.5, puts the flavored quiver image inside
that pure-gauge algebra, using dressed-minuscule generation and the formulas
(2.7)--(2.8). Only this inclusion is used, not the converse characterization.

The quiver in~\eqref{eq:datum} is loop-free for every $m\geq1$, including the
unframed case $m=1$. Dressed-minuscule generation and the indicated inclusion do
not require quiver balance. Apply these statements while keeping the source
characters and the zero first-node shift. The source convention satisfies
$D_\mu w_\lambda=q^{2\langle\mu,\lambda\rangle}w_\lambda D_\mu$; hence $q=u$,
not $u^2$. This is exactly~\eqref{eq:operators}. At the wall $b/a=u^{2h}$, its
reflection-shift partner is
\begin{equation}\label{eq:partner}
 (R,S)\longmapsto(S+h,R-h).
\end{equation}
The source uses the internal zero-section convention, which is the convention reached
in~\eqref{eq:rawlower}. The simultaneous target negation preserves $b/a$, commutes
with Weyl reflection and every shift, and is an automorphism of the integral Laurent
ring. The scalar $-u^{1-k}/x$ is common to the entire fixed class, as is the optional
component unit $\tau^{N-2}$. They preserve the paired residues. Tensoring by the
electric line was performed geometrically before applying the image theorem; no
coefficientwise arbitrary renormalization is used here.

It follows that the assembled coefficients $B_R$ have only simple even-root poles
and satisfy~\eqref{eq:mainresidue} for all $N,k$. The fixed lattice labels of the
reduced closure are exactly $R=-1,\ldots,N-1$ with second entry $N-2-R$, which
explains the $N+1$ terms.

\begin{corollary}\label{cor:poles}
Up to Laurent units, a squarefree denominator valid for $B_R$ is
\begin{equation}\label{eq:denombound}
 D_R=\prod_{\substack{R-N+1\leq h\leq R+1\\h\ne2R-N+2}}
                  (1-u^{-2h}b/a).
\end{equation}
It has $N$ factors; further cancellations may occur.
\end{corollary}
\begin{proof}
The partner index is $R'=N-2-R+h$. A pole can remain only if $R'$ is in the support
and $R'\ne R$. An absent partner forces zero residue. A self-pair gives twice the
residue, hence zero in characteristic zero. Since all poles are at most simple,
these vanishing residues remove the corresponding walls. The inequalities
$-1\leq R'\leq N-1$ give the range in~\eqref{eq:denombound}, and the omitted value
is exactly the self-wall. This includes $N=2$ and the even-$N$ central label.
\end{proof}

\subsubsection{Divisibility over the integral ring}
Put $\mathcal A=\Z[u^{\pm1},x_1^{\pm1},\ldots,x_m^{\pm1}]$ and
$L=\mathcal A[a^{\pm1},b^{\pm1}]$. The integral image theorem and the explicit
descent in Lemma~\ref{lem:degree} put each numerator in $L$, with only root-binomial
denominators and no inverted nonunit integers. The binomials
$p_h=1-u^{-2h}b/a$ are primitive prime elements: $L/(p_h)$ is the Laurent domain
obtained by $b=u^{2h}a$. For formal $u$, distinct $h$ give nonassociate primes.
Simplicity of poles therefore permits a common squarefree denominator $D$.

Take $F\in L^{\mathfrak S_2}$. At a wall $p_h=0$, the shifted arguments for $R$ and
$R'=N-2-R+h$ are exchanged. Indeed their unordered pairs are both
\[
 \{u^{2R}a,\ u^{2(N-2-R+h)}a\}.
\]
Symmetry of $F$ makes the paired values equal. Multiplication by $p_h$ followed by
restriction to the wall gives zero by~\eqref{eq:mainresidue}; absent and self-paired
terms vanish separately.

Write $\T_{N,k}F=P/D$, with $P\in L$. The preceding vanishing says that $P$ maps
to zero in the fraction field of the domain $L/(p_h)$, hence in that domain itself.
Thus $p_h$ divides $P$ in $L$, not only after tensoring with $\C$ or $\Q$. Every prime
factor of the squarefree $D$ divides $P$, so unique factorization gives $D\mid P$.
Consequently $\T_{N,k}F\in L$. Weyl covariance of the coefficients and shifts gives
$\T_{N,k}F\in L^{\mathfrak S_2}$. This proves the residue and integral-action assertions
of Theorem~\ref{thm:main}.

The source specialization may now be made in this polynomial answer and preserves
integrality. No assertion is made about separately evaluating singular rational
summands at roots of unity. Likewise, an electric recurrence is consistent with but
not needed for this proof. For a single matter exponent $m=k-r$, increasing $m$ by
$N$ multiplies the ordinary coefficient at $R$ by
\begin{equation}\label{eq:recurrence}
 u^{-2R^2+2(N-2)R+N-2}a^{N-1-R}b^{R+1}.
\end{equation}
The paired multipliers agree at $b/a=u^{2h}$ because their exponent difference is
$(R-R')[-2(R+R')+2(N-2)+2h]=0$. This propagates an established residue identity in
both electric directions; it is not a replacement for the uniform image-theorem
argument or a claim of electric periodicity.

\subsection{Triangular characterization and numerator structure}\label{sec:structural}
The even-root image theorem also yields a useful uniqueness statement on precisely
these supports. Its domain is essential: permitting extra odd-root poles in arbitrary
rational coefficients would not give the following argument.

\begin{lemma}\label{lem:maximalsupport}
Let $\Delta$ be a finite Weyl-covariant difference operator in a fixed determinant
component. Suppose its coefficients lie in $L$ localized only at even-root binomials,
have at most simple poles, and satisfy the paired-residue identities. If $(R,S)$,
$R\geq S\geq0$, is a maximal dominant label in its nonzero support, then its symmetric
coefficient $f_{R,S}$ belongs to $L$.
\end{lemma}
\begin{proof}
Maximality and Weyl covariance place every first exponent in the interval $[S,R]$.
At $p_h=0$, the partner of $(R,S)$ has first exponent $S+h$. If $h<0$ or $h>R-S$,
it is absent, and the residue vanishes. If $h=R-S$, it is self-paired, so twice the
residue vanishes and hence the residue is zero. Simplicity removes these poles.
The only remaining denominator primes divide $V_{R-S}=\prod_{h=0}^{R-S-1}p_h$.
Primitivity and unique factorization give $V_{R-S}B_{R,S}\in L$ integrally.
Multiply by $(-u)^S$ and undo the shift $(a,b)\mapsto(u^Ra,u^Sb)$ in
\eqref{eq:operators}; these operations preserve $L$ and give $f_{R,S}\in L$.
The same proof includes $R=S$, where $V_0=1$.
\end{proof}

\begin{corollary}\label{cor:unique}
In the actual integral Coulomb algebra with zero first-node magnetic charge and
support contained in $\overline{\cO_{(N-1,-1)}}$, there is at most one normalized
operator with prescribed top symmetric coefficient and its Weyl mate whose lower
symmetric coefficients are bar-fixed and $O(u)$. More generally, on this first-strip support, uniqueness holds in
the explicitly restricted integral even-root, simple-pole residue class.
\end{corollary}
\begin{proof}
Subtract two candidates. Their top orbit cancels. A nonzero maximal dominant label
of the difference therefore has $R\geq S\geq0$. Lemma~\ref{lem:maximalsupport}
makes its symmetric coefficient a finite Laurent polynomial in $u$, with coefficients
in the Laurent ring of the other characters. Positive order removes all powers
$u^d$ for $d\leq0$. Bar invariance identifies the coefficients of $u^d$ and $u^{-d}$,
so every remaining coefficient vanishes too, a contradiction. The abelianization is
injective, so equality of operators is equality of the algebra classes.
\end{proof}

Existence is supplied by~\eqref{eq:mainformula}, whose $O(u^j)$ bound is stronger
than required. Thus the normalized geometric classes have an intrinsic triangular
characterization throughout the fixed infinite family. The common component unit is
removed uniformly before applying the integral-$u$ argument. This is not uniqueness
in a larger rational class, nor an exhaustion of a canonical basis.

\begin{corollary}\label{cor:positive}
Put $\rho=b/a$ and define
\begin{equation}\label{eq:positive}
 \begin{aligned}
 \Delta_j(u,\rho)&=\prod_{l=1}^j
       (1-u^{N-2j+2l}\rho)(1-u^{N-2j+2l}\rho^{-1}),\\
 N_i(u,\rho)&=\Delta_j(u,\rho)G_i(u,\rho),\\
 \mathcal P_{j,k}(u,a,b;x_r)&=
       \sum_{r=0}^me_r\delta_j^{k-r}N_{r-k}(u,b/a).
 \end{aligned}
\end{equation}
Then $f_{\lambda_j}^{(k)}=(-1)^ju^j\mathcal P_{j,k}/(x\Delta_j)$, and
$\mathcal P_{j,k}$ has nonnegative integer coefficients in $u$ and Laurent source
and target characters. With all characters other than $u$ fixed,
\begin{equation}\label{eq:palindromic}
 \mathcal P_{j,k}(u^{-1})=u^{-2j(N-j)}\mathcal P_{j,k}(u).
\end{equation}
Its scalar degrees lie in $[0,2j(N-j)]$, and its finite total multiplicity is
\begin{equation}\label{eq:mass}
 \left.\mathcal P_{j,k}\right|_{u=a=b=1}
       =\binom Nj\prod_{r=1}^m(1+x_r).
\end{equation}
\end{corollary}
\begin{proof}
By~\eqref{eq:free}, $N_i$ is the character of a homogeneous free basis over the
off-diagonal parameter ring. Its coefficients are nonnegative integers and its
scalar degrees are nonnegative. The deck-compensation calculation in
Lemma~\ref{lem:degree} cancels all fractional target powers term by term; the
$e_r$ have nonnegative coefficients as well. Since
$\Delta_j(u^{-1},\rho)=u^{-2j(N-j+1)}\Delta_j(u,\rho)$,
the fixed-sector identity in~\eqref{eq:threeidentities} gives~\eqref{eq:palindromic}.
Nonnegative degree and reciprocity give the upper bound. Every $N_i$ has total
multiplicity $\binom Nj$ by Theorem~\ref{thm:shuffle}; setting the finite numerator
characters to one gives~\eqref{eq:mass}.
\end{proof}

The stated interval is a bound, not a claim that either endpoint occurs. Positivity
refers to this specified common homogeneous-parameter denominator and extracted
global sign; it is not asserted after arbitrary denominator reduction. The shuffle
formula and freeness do not by themselves supply a cancellation-free combinatorial
basis or a positive multiplication rule.

Finally the complete electric family can be reconstructed from $N$ consecutive
integer labels. Directly from~\eqref{eq:mainformula}, including the separately
computed top term, one has
\begin{equation}\label{eq:reconstruct}
 f_{\lambda_j}^{(k+N)}=a^jb^{N-j}f_{\lambda_j}^{(k)},
 \qquad 0\leq j\leq\lfloor N/2\rfloor.
\end{equation}
The multipliers are units, so reconstruction works in both directions. In ordinary
coordinates the multiplier is~\eqref{eq:recurrence}; it depends on the magnetic
label. This is electric reconstruction, not periodicity and not one common scalar
multiple of the whole operator. These arguments complete the proof of Theorem~\ref{thm:main}.

\section{Constant-sign cones and first-strip examples}\label{sec:cone}\label{sec:same-sign}
In this section the number of matter columns is four, and the first-node magnetic charge is zero.  Central translation identifies the magnetic bases and their electric covers, but does not preserve the matter incidence.  Recomputing that incidence nevertheless gives a uniform answer on both constant-sign cones.

For a dominant pair $R\geq S$, put
\begin{equation}\label{cone:normalizers}
\begin{split}
 V_{R-S}(A,B)&=\prod_{i=0}^{R-S-1}(1-u^{-2i}B/A),\\
 P_d(Z)&=\prod_{v=1}^4\prod_{i=1}^d(1+u^{-(2i-1)}Z/x_v),\qquad P_0=1,\\
 J_{R,S}(A,B)&=(-u)^{-S}
 \frac{P_{\max(-R,0)}(A)P_{\max(-S,0)}(B)}{V_{R-S}(A,B)}.
\end{split}
\end{equation}
Other labels are obtained by Weyl reflection.  Formula \eqref{cone:normalizers} is the extension of \eqref{eq:J} to arbitrary signs: the even-root factors and negative odd matter factors in \cite[eq.~(53)]{Gaiotto} give the displayed products after removal of the common determinant-component factor.  It agrees with the positive-cone normalization and with $J_{1,-1}=-uP_1(B)/V_2$; at the negative central label it gives $J_{-1,-1}=-uP_1(A)P_1(B)$.

Fix $r>s$ with either $s\geq0$ or $r\leq0$, and write
\[
 \nu=r-s,\quad c=r+s,\quad \lambda_j=(r-j,s+j),
 \qquad 0\leq j\leq\lfloor\nu/2\rfloor.
\]
Let $H_i^{\nu,j}$ denote the saturated electric-cover character of Section~\ref{sec:characters} for length $\nu$, and set
\begin{equation}\label{cone:characters}
 \delta_j^\nu=a^jb^{\nu-j},\qquad
 G_i^{\nu,j}(u,b/a)=H_i^{\nu,j}(u^2,u^{\nu-2j}b/a).
\end{equation}
Only the index of $H_i$ is read modulo $\nu$.  Let $x=\prod_{v=1}^4x_v$, and let $\chi$ be an arbitrary Laurent character of the source torus.

\begin{theorem}\label{cone:theorem}
For $H,K\in\Z$, there is a rephased Cautis--Williams simple on the closure of $\lambda_0$ with symmetric coefficients
\begin{equation}\label{cone:formula}
\begin{split}
 f_{r,s}&=a^Hb^{K+H}\chi,\\
 f_{r-j,s+j}&=(-u)^j(ab)^H\chi\,\delta_j^K
                 G_{-K}^{\nu,j}(u,b/a)\qquad(j>0).
\end{split}
\end{equation}
The geometric electric and source labels producing these coefficients are
\begin{equation}\label{cone:relabel}
 \begin{cases}
 \ell^K(\det V_0)^H\chi,&s\geq0,\\
 \ell^{K+4\nu}(\det V_0)^{H-4r}x^{c}\chi,&r\leq0.
 \end{cases}
\end{equation}
Here $\ell$ is the leading-image line of the translated lattice, with upper character $b$.  These labels exhaust all electric characters of the noncentral reductive stabilizer.  The lower coefficients are bar-fixed and $O(u^j)$, and they descend to integral target characters.  They are the unique support-bounded triangular lifts with these properties inside the actual integral Coulomb algebra.  A central lower label is counted once.
\end{theorem}

\begin{proof}
If $s\geq0$, all lattices satisfy $L\subset\cO^2$.  The external image $\cO^2/L$ is a vector bundle of rank $c$; the internal incidence has no constraints.  Tensor the saturated electric-cover extension with its four-column regular Koszul complex.  All exterior terms are retained categorically.  On a transverse slice of dimension $2j$, its coefficients are tensor products of the established MCM electric modules with vector bundles, hence are MCM.  The bounded coherent zero-matter detector of Proposition~\ref{prop:simple} therefore gives the strict restriction and costalk inequalities, and proves simplicity.  The Euler factors cancel only when passing from external to internal polarization.

If $r\leq0$, all lattices contain $\cO^2$.  The external image is zero.  The incidence simple is the regular polynomial matter bundle with the excess dg generators acting by zero, not the exterior structure sheaf of the excess derived intersection.  Tensoring the electric extension with this polynomial bundle again preserves every transverse MCM bound, so the same detector proves simplicity.  For a geometric label $\ell^k(\det V_0)^h$ the prescribed shift is $[\nu/2-k-2h]$; its difference from the component baseline $c/2$ is the integer $-s-k-2h$.

Here is the normalization calculation, including the negative matter factors.  At $\lambda_j$, the translated quotient has chains of lengths $\nu-j,j$, and its cyclic determinant gives
\[
 c_j^\nu=A^jB^{\nu-j}u^{2j(\nu-j)},\qquad
 \operatorname{ch}(\ell^{[k]})=c_j^kH_{-k}^{\nu,j}(u^2,B/A).
\]
A common scalar is $\kappa=(-u)^{-s}u^{sk+ch}$.  With $A=u^{-r+j}a$, $B=u^{-s-j}b$, one has
\[
 c_j=u^{-s}\delta_j,\qquad AB=u^{-c}ab.
\]
For the nonnegative cone the ordinary coefficient is
$\kappa(AB)^h\chi c_j^kH_{-k}/V_{\nu-2j}$.
Dividing by $J_{r-j,s+j}$ gives \eqref{cone:formula} with $k=K,h=H$.

For $R,S\leq0$, the internal constrained modes instead contribute
\[
 E_{R,S}=\prod_{v=1}^4
 \left[\prod_{i=1}^{-R}(1+u^{2i-1}x_v/A)
       \prod_{i=1}^{-S}(1+u^{2i-1}x_v/B)\right].
\]
Indeed a constrained mode numbered $q$, with $0\leq q<-R$, has weight
\[ u^{-1-2R-2q}x_v/A. \]
This gives the first product by $i=-R-q$.
Factoring each $1+y$ against $1+y^{-1}$ yields the exact identity
\begin{equation}\label{cone:euler-ratio}
 \frac{E_{R,S}}{P_{-R}(A)P_{-S}(B)}
 =u^{4(R^2+S^2)}x^{-R-S}A^{4R}B^{4S}.
\end{equation}
All scalar powers of $u$ on the right cancel under the symmetric substitution.  Thus the ordinary coefficient
\[
 \frac{\kappa x^c\chi E_{r-j,s+j}(AB)^hc_j^k
 H_{-k}^{\nu,j}(u^2,B/A)}{V_{\nu-2j}}
\]
has symmetric coefficient
$(-u)^j\chi(ab)^ha^{4(r-j)}b^{4(s+j)}\delta_j^kG_{-k}$.
Putting $k=K+4\nu$, $h=H-4r$ gives \eqref{cone:formula}, since
\[
 \delta_j^{4\nu}=a^{4j}b^{4(\nu-j)},\qquad
 (ab)^{-4r}a^{4(r-j)}b^{4(s+j)}\delta_j^{4\nu}=1.
\]
This also explains why the negative-cone electric relabeling cannot be omitted.  The parity of the prescribed shift and simultaneous target-character negation gives the same component sign as $\kappa$; no numerical half-integral sign is required.

The balanced cover-character theorem now gives bar invariance, order $j$, integral descent, and, after extracting $(-1)^j$, nonnegative integral numerators over its homogeneous-parameter denominator.  Algebra membership comes from the constructed simple; its residue conditions follow from the image theorem of \cite{SS}.

For uniqueness, let a nonzero difference of two such actual algebra classes have maximal nonzero magnetic stratum $(R,S)$.  Closed/open localization in equivariant $K$-theory removes the larger strata on which its restriction vanishes.  On this orbit the classical incidence is a vector bundle, so its restriction is a finite Laurent stabilizer character $\psi$.  Its ordinary coefficient is $\psi/V_{R-S}$ in the nonnegative cone and $E_{R,S}\psi/V_{R-S}$ in the nonpositive cone, up to character units.  Formula \eqref{cone:euler-ratio} shows that the corresponding symmetric coefficient is Laurent in either case.  A finite Laurent polynomial fixed by $u\mapsto u^{-1}$ and $O(u)$ is zero, a contradiction.  Induction over the finite support proves uniqueness.  In the negative cone this argument uses actual algebra membership, not merely an even-root residue ansatz.
\end{proof}

At a central charge $(s,s)$ the reductive target stabilizer is $\GL_2$, not a torus.  Its complete electric labels are
\[
 V_{q,H}=\det^H\Sym^qV_0,\qquad q\geq0,\ H\in\Z,
 \qquad \chi_{q,H}(a,b)=(ab)^H\sum_{i=0}^q a^{q-i}b^i.
\]
For $s\geq0$, use $V_{q,H}\chi$; for $s<0$, use
$\det^{H-4s}\Sym^qV_0\otimes x^{2s}\chi$.  If $w$ is the central target weight of this geometric representation, the categorical shift is $[-w]$ and the common scalar is $(-u)^{-s}u^{sw}$.  The same Euler ratio proves that the sole symmetric coefficient is $\chi_{q,H}(a,b)\chi$.  These are simple because their open stabilizer representations are irreducible and there is no magnetic boundary.  The representations extend by evaluation of the target jet group at zero.  No cyclic cover of length zero is used, and no analogous simplicity claim is made for a reducible torus restriction at a noncentral orbit.

\subsection{The first bubble without fractional characters}\label{subsec:first-bubble}
We return to arbitrary $m\geq1$ and the first-strip support of
Theorem~\ref{thm:main}. The first degeneration has transverse surface
\[
 X_1=\Spec\C[\alpha,c,d]/(cd+\alpha^N),\qquad N=n+1.
\]
This is the type-$A_{N-1}$ minimal-degeneration slice~\cite{MOV}. Its actual
cyclic normalization cover has coordinates $p,q$, with
$\alpha=pq$, $c=p^N$, $d=-q^N$ and deck action
$(p,q)\mapsto(\zeta p,\zeta^{-1}q)$. In centered variables choose
$\xi^N=a/b$; the function weights on the cover are $u\xi,u\xi^{-1}$.
This supplies a direct small-dimensional realization of the normalization used
in the general shuffle formula.

For $s\in\Z/N\Z$ put
\begin{equation}\label{example:finite-residue}
 S_s(u,\xi)=
 \sum_{\substack{0\leq\alpha,\beta<N\\\alpha-\beta\equiv s\pmod N}}
       u^{\alpha+\beta}\xi^{\alpha-\beta},\qquad
 D_N=(1-u^Na/b)(1-u^Nb/a).
\end{equation}
Writing each monomial exponent modulo $N$ shows that the sector character is
$G_s=S_s/D_N$. Since $\delta_1=b\xi$, the general electric coefficient is
\begin{equation}\label{example:first-all-k}
 f_{N-2,0}^{(k)}
 =-\frac{u}{xD_N}\sum_{r=0}^m e_r b^{k-r}\xi^{k-r}S_{r-k}(u,\xi).
\end{equation}
The full integer $k$ remains in the electric power.

\begin{corollary}\label{example:root-free}
At electric-zero normalization $k=m$, the first bubbling coefficient equals
\begin{equation}\label{example:root-free-formula}
 -\frac{u}{xD_N}
 \sum_{r=0}^m e_r
 \sum_{\substack{0\leq\alpha,\beta<N\\
                 \alpha-\beta\equiv r-m\pmod N}}
 u^{\alpha+\beta}a^{\kappa}b^{m-r-\kappa},
 \qquad \kappa=\frac{m-r+\alpha-\beta}{N}.
\end{equation}
Its numerator after extraction of $-u/(xD_N)$ is a polynomial, with
nonnegative integer coefficients, in $a,b,u$ and the elementary source characters.
\end{corollary}
\begin{proof}
The residue condition makes $\kappa$ integral. As $0\leq\alpha,\beta<N$,
\[
 -N<m-r+\alpha-\beta<(m-r+1)N.
\]
It follows that $0\leq\kappa\leq m-r$. Thus the target exponents in
\eqref{example:root-free-formula} are both nonnegative and all fractional
characters cancel term by term. For arbitrary $k$, the same calculation proves
integral Laurent descent, but not nonnegative target exponents.
\end{proof}

For $m=4,N=2$ this reads
\begin{equation}\label{example:small-first}
 -\frac{u}{xD_2}\bigl\{
 (1+u^2)(a^2b^2+e_2ab+e_4)
 +u(a+b)(e_1ab+e_3)\bigr\}.
\end{equation}
For $N=3$ the full operator has four labels; for $N=4$ it has five. Their residue
cancellations are consequences of the full-support theorem. The first-lower
formula alone would not justify a polynomial-preserving operator when deeper
magnetic strata are present.

\subsection{Composite deck order and the central slice}\label{subsec:small-characters}
The normalization formula has no prime-$N$ restriction. At $N=4$, the scalar
shifts are $h_1=h_3=3/2$ but $h_2=1$. At the central slice $j=2$, write
$\rho=s^2$ and
\[
 \mathcal D=(1-ts^2)(1-ts^{-2})(1-t^2s^2)(1-t^2s^{-2}).
\]
The six descending shuffles give
\begin{equation}\label{example:four-sectors}
 \begin{aligned}
 H_0&=(1+t+2t^2+t^3+t^4)/\mathcal D,\\
 H_1=H_3&=t(1+t+t^2)(s+s^{-1})/\mathcal D,\\
 H_2&=\bigl(t+2t^2+t^3+t^2(s^2+s^{-2})\bigr)/\mathcal D.
 \end{aligned}
\end{equation}
Every numerator has total multiplicity $\binom42=6$. These are specializations
of Theorem~\ref{thm:shuffle}, not input for its proof. At $N=6$, the nonzero-sector
shifts are $5/2,2,3/2,2,5/2$; the two deepest ranks are
$\binom62=15$ and $\binom63=20$. The greatest-common-divisor correction is already
necessary in these small examples.

The same length-four magnetic base occurs after central translation in the next
section. Its matter map changes, however:
\begin{center}
\begin{tabular}{@{}lccc@{}}
\toprule
Matter map & Open stratum & Middle stratum & Central stratum\\
\midrule
$\im T$ for $(3,-1)$ & $3$ & $2$ & $2$\\
$\im T^2$ for $(2,-2)$ & $2$ & $1$ & $0$\\
\bottomrule
\end{tabular}
\end{center}
These are the ranks of the actual images on Jordan types $(4)$, $(3,1)$, and
$(2,2)$. The electric-cover characters remain available, but they cannot determine
the saturated exterior powers of the new image by themselves.

\section{The first mixed-depth boundary correction}\label{sec:boundary}
We retain four matter columns and first-node charge zero, but now use the original regular matter lattice on
\[
 \overline{\mathcal O}_{(2,-2)}
 =\mathcal O_{(2,-2)}\sqcup\mathcal O_{(1,-1)}\sqcup\mathcal O_{(0,0)}.
\]
The matter image is $\im T^2$, not $\im T$.  Its saturated exterior coefficients need not be MCM.  We construct the resulting boundary object with its actual dg actions before computing its localization.

\subsection{The image bundle and its finite cover}
On the central chart write $g=1+\mathsf A/z+\mathsf B/z^2$.  The equations are
\begin{equation}\label{boundary:chart}
 \operatorname{tr}\mathsf A=0,\quad
 \operatorname{tr}\mathsf B+\det\mathsf A=0,\quad
 \operatorname{tr}(\mathsf A\mathsf B)=0,\quad \det\mathsf B=0.
\end{equation}
They define a normal four-dimensional complete intersection $X$.  Indeed, after eliminating $\operatorname{tr}\mathsf A$, there are three equations in seven variables.  For $\mathsf A$ with eigenvalues $\pm\alpha$, $\alpha\ne0$, the fiber is
$\mathsf B_{11}=\mathsf B_{22}=\alpha^2/2$,
$\mathsf B_{12}\mathsf B_{21}=\alpha^4/4$.
The fibers over nonzero nilpotent $\mathsf A$ and over zero have dimensions one and two, respectively.  Thus there is a unique generic component and no additional component; complete-intersection unmixedness gives reducedness.  The complement of the smooth open orbit has dimension two, so Serre's criterion gives normality.  The three orbit loci are $\mathsf B\ne0$, $\mathsf B=0\ne\mathsf A$, and the origin.

For an original matter column $s=s_0+zs_1$, subtracting
$g(s_0-\mathsf A s_1)+zg s_1$ gives the genuine quotient equations
\begin{equation}\label{boundary:mattermatrix}
 \binom{-\mathsf A s_0+\operatorname{adj}(\mathsf B)s_1}
 {-\mathsf B s_0+\mathsf B\mathsf A s_1}
 =\mathsf M\binom{s_0}{s_1},\qquad
 \mathsf M=\begin{pmatrix}-\mathsf A&\operatorname{adj}(\mathsf B)\\
                         -\mathsf B&\mathsf B\mathsf A\end{pmatrix}.
\end{equation}
Its ranks are $2,1,0$ on the three strata.  There are sixteen ambient odd generators, but only eight intrinsic ones on the upper orbit.  In particular, the first-strip free-generator splitting does not apply.

Here is a finite cover of the magnetic chart, with the linearizations needed below.  Set
\[
 \begin{gathered}
 \epsilon(p)=(-p_2,p_1),\quad q^\sharp=(q_2,-q_1)^t,\\
 \mathbf z=(p_1,p_2,\rho)^t,\quad \mathbf w=(q_1,q_2,\sigma),
 \qquad qp+\rho\sigma=0.
 \end{gathered}
\]
The rank-one traceless matrix $\mathbf z\mathbf w$ parametrizes the minimal nilpotent cone $Y\subset\mathfrak{sl}_3$; the pairs are taken modulo
$(\mathbf z,\mathbf w)\mapsto(t\mathbf z,t^{-1}\mathbf w)$.
Put
\begin{equation}\label{boundary:covercoordinates}
 \begin{gathered}
 \zeta=-\rho\sigma/2,\quad
 \mathsf A=pq+\rho\sigma I_2/2,\\
 v=\sigma p/2-\rho q^\sharp/4,\quad
 w=-\rho q/4-\sigma\epsilon(p)/2,\quad \mathsf B=vw.
 \end{gathered}
\end{equation}
The generator
\begin{equation}\label{boundary:deck}
 \gamma:(p,\rho;q,\sigma)\longmapsto
 (q^\sharp,-2i\sigma;\epsilon(p),-i\rho/2)
\end{equation}
induces $(\mathsf A,v,w,\zeta)\mapsto(\mathsf A,iv,-iw,-\zeta)$.
It realizes $X=Y/\mu_4$, and is free over $U=\{\mathsf B\ne0\}$.
For completeness, the rank-one equations in the matrix
\[
 \begin{pmatrix}\mathsf A+\zeta I_2&v-w^\sharp\\
             -2(\epsilon(v)+w)&-2\zeta\end{pmatrix}
\]
give
$\det\mathsf A=-\zeta^2$, $wv=\zeta^2$,
$\mathsf Av=-\zeta w^\sharp$, $w\mathsf A=\zeta\epsilon(v)$,
and $\zeta\mathsf A=v\epsilon(v)-w^\sharp w$.
These reduce every invariant monomial to entries of $\mathsf A,vw$: mixed monomials factor through $vw$, terms with one $\zeta$ reduce by the two middle identities, and pure invariant quartics reduce by the last identity.  Hence the invariant ring is generated by $\mathsf A,\mathsf B$.  Both it and the domain defined by \eqref{boundary:chart} have dimension four, so the resulting surjection is an isomorphism.  The cone $Y$ is normal and CM: it is the trace hyperplane in the rank-one matrix cone, whose diagonal section ring is CM~\cite{BV} by the K\"unneth calculation on $\mathbb P^2\times\mathbb P^2$.  It is smooth off its vertex.  The displayed identities also show that only the vertex lies over the origin and that every nontrivial deck stabilizer lies over $\mathsf B=0$.

The crucial image calculation is the factorization
\begin{equation}\label{boundary:factorization}
 \mathsf M=FG,\qquad
 F=\begin{pmatrix}p&q^\sharp\\-\rho v/2&-\sigma v\end{pmatrix},\qquad
 G=\begin{pmatrix}-q/2&(\sigma/2)\epsilon(v)\\
                 -\epsilon(p)/2&(\rho/4)\epsilon(v)\end{pmatrix}.
\end{equation}
Multiplication verifies the identity modulo $qp+\rho\sigma$.
The minors of $F$ include nonzero scalar multiples of
$\rho\sigma,v_1^2,v_1v_2,v_2^2$; thus $F$ has rank two over the preimage of $U$ and its image equals that of $\mathsf M$.  Its columns have gauge weights $+1,-1$ and are exchanged exactly by $\gamma$.

Let
$\mathcal F=\{([\mathbf z],[\mathbf w]):\mathbf w\mathbf z=0\}
\subset\mathbb P^2\times\mathbb P^2$.
The punctured cone $Y^*$ is the punctured total space of
$\mathcal O_{\mathcal F}(-1,-1)$.
Consequently the saturated image pulls back to the sum of the two opposite gauge line bundles $\mathcal L,\mathcal L^{-1}$ on $Y^*$.
The complement of the preimage of $U$ has codimension two, so this identifies the reflexive extensions, not merely the generic fibers.
The electric frame $v$ has deck weight $+1$; its $k$th power selects sector $-k\bmod4$ after descent.

Let $\mathcal M_k$ denote the ordinary termwise saturated intrinsic image-Koszul module with electric/source twist $\ell^k/x$, where $x=\prod_{j=1}^4x_j$.  Explicitly, extend each open exterior coefficient reflexively and extend the ambient wedge actions and contractions uniquely.  Their identities hold on the dense open set and hence globally in the torsion-free target modules.  These extensions glue equivariantly, and produce a bounded coherent ambient dg module.  Write $T_l$ for its exterior-degree-$l$ coefficient, with all electric and source twists retained.

\subsection{The precise depth defect}
For $m\geq0$, the line module associated with $\mathcal O_{\mathcal F}(m,0)$ is
\[
 \mathscr Q_m=\bigoplus_{n\in\Z}
 H^0(\mathcal F,\mathcal O_{\mathcal F}(m+n,n)).
\]
The opposite line is obtained by interchanging the factors.  Its intermediate vertex local cohomology is computed by the incidence sequence
\begin{equation}\label{boundary:incidence}
 0\longrightarrow\mathcal O(m+n-1,n-1)
 \longrightarrow\mathcal O(m+n,n)
 \longrightarrow\mathcal O_{\mathcal F}(m+n,n)\longrightarrow0
\end{equation}
on $\mathbb P^2\times\mathbb P^2$.  The complete list for $0\leq m\leq4$ is
\begin{equation}\label{boundary:bwbtable}
\begin{array}{c|c|c}
 m&H^1(\mathcal F,\mathcal O(m+n,n))&H^2(\mathcal F,\mathcal O(m+n,n))\\ \hline
 0,1,2&0\text{ for all }n&0\text{ for all }n\\
 3&\C\text{ at }n=-2&\C\text{ at }n=-3\\
 4&\C^3\text{ at }n=-2&\C^3\text{ at }n=-4
\end{array}
\end{equation}
All omitted groups vanish.  This is also the relevant Borel--Weil--Bott calculation, but the elementary sequence fixes its equivariance: for $m=4,n=-2$ it identifies $H^1$ with
$H^0(\mathbb P^2,\mathcal O(1))\otimes H^2(\mathbb P^2,\mathcal O(-3))$.
At $m=4,n=-3$, the only possible map is the isomorphism
$H^2(\mathcal O(0,-4))\to H^2(\mathcal O(1,-3))$ induced by contraction, so there is no missing intermediate group.

Before descent, exterior degree $l$ decomposes into opposite-line summands indexed by
$0\leq a,b\leq4$, $a+b=l$, with multiplicity
$\binom4a\binom4b$ and imbalance $a-b$.
Thus $H^2$ at the vertex occurs only in degrees $3,4,5$.
The square of \eqref{boundary:deck} acts on both triples by
$D_0=\operatorname{diag}(1,1,-1)$, with no additional gauge sign.
The top-cohomology factor $H^2(\mathcal O(-3))$ has character
$\det D_0=-1$.  Hence its tensor with $H^0(\mathcal O(1))$ has square-eigenvalues $(-1,-1,+1)$.  Interchanging the two extreme summands gives the following deck multiplicities:
\begin{equation}\label{boundary:decktable}
\begin{array}{c|c|rrrr}
 l&(a,b)&1&i&-1&-i\\ \hline
 3&(3,0),(0,3)&0&4&0&4\\
 4&(4,0),(0,4)&1&2&1&2\\
 5&(4,1),(1,4)&0&4&0&4
\end{array}
\end{equation}
For degrees three and five the source in \eqref{boundary:incidence} is $\mathcal O(0,-3)$; its square-character is $-1$, explaining the odd sectors.  The same ungraded deck table occurs in local $H^3$, at the different diagonal degrees of \eqref{boundary:bwbtable}; it lies one costalk degree higher.

\begin{proposition}\label{boundary:depthfailure}
For every $k\in\Z$, the saturated coefficient $T_4$ has depth exactly two at the central point.  For odd $k$, $T_3$ and $T_5$ also have depth two.  Nevertheless $\mathcal M_k[2-k]$ lies in the bounded coherent Koszul-perverse heart.  It has no boundary-supported quotient and no middle-stratum subobject.
\end{proposition}
\begin{proof}
Reflexivity gives depth at least two.  Local cohomology commutes with the finite pushforward and characteristic-zero deck projection, and every sector in the middle row of \eqref{boundary:decktable} is nonzero.  This proves the depth assertion for every integer electric label, not only for one congruence class.  All local $H^2$ groups are concentrated at a single diagonal degree and are annihilated by the magnetic maximal ideal.

Away from the vertex, the line bundles on the smooth cover give MCM modules after finite descent.  In particular the transverse middle-stratum depth is two.  Under the bounded coherent zero-matter detector, the coefficient $T_l[l+2-k]$, of regrading weight $k-l$, becomes $T_l[2]$.  At the center its restriction is in degrees at most $-2$ and its costalk in degrees at least zero.  These are the heart inequalities, with a strict restriction inequality.  At the middle stratum the transverse bounds are both strict.  Exactness and conservativity of the detector establish the assertions.  Depth two is enough for heart membership here; it does not give the strict central costalk bound required for simplicity.
\end{proof}

\subsection{The actual boundary row and its canonical lift}
Let $C_l=H^2_{\{0\}}(T_l)$.  Tensoring these finite-dimensional magnetic modules with the original matter polynomial ring gives a finite coherent row with its inherited ambient action and differential.  These maps must not be replaced by independent point-sheaf injections.

Write $\alpha_j,\beta_j$ for the two intrinsic dual-image generators in column $j$.  Formula \eqref{boundary:factorization} gives the ambient wedge operations
\begin{equation}\label{boundary:rowactions}
 p_i\alpha_j+q_i^\sharp\beta_j,\qquad
 -\rho v_i\alpha_j/2-\sigma v_i\beta_j,
\end{equation}
and the matter contractions
\begin{equation}\label{boundary:rowdifferential}
 d\alpha_j=-q s_{0,j}/2+(\sigma/2)\epsilon(v)s_{1,j},\qquad
 d\beta_j=-\epsilon(p)s_{0,j}/2+(\rho/4)\epsilon(v)s_{1,j}.
\end{equation}
Apply the electric deck projector to these formulas, with the usual wedge signs.  For odd $k$, the row has dimensions $4,2,4$ and its nonzero maps have shape
\begin{equation}\label{boundary:rowshape}
 C_3\xrightarrow{\text{ambient exterior action}}C_4
 \xleftarrow{\text{matter differential}}C_5.
\end{equation}
Indeed multiplication by $p_i$ on
$H^1(\mathcal O_{\mathcal F}(1,-2))$ is multiplication on
\[ H^0(\mathcal O(0))\otimes H^2(\mathcal O(-3)) \]
in
\eqref{boundary:incidence}.  The two top coordinates give the nonzero odd-sector maps into
\[ H^1(\mathcal O_{\mathcal F}(2,-2)). \]
All exterior arrows out of $C_4$ land in an acyclic line bundle, as do the bottom-row cubic terms.  The only surviving contractions, on the positive and negative extreme degree-five pieces, are respectively
\[
 -\epsilon(p)s_{0,j}\alpha_1\alpha_2\alpha_3\alpha_4/2,
 \qquad -q s_{0,j}\beta_1\beta_2\beta_3\beta_4/2.
\]
Every other contraction has imbalance at most two or an acyclic diagonal degree.  Thus every path of length two vanishes.  For even $k$ the row consists of a single one-dimensional $C_4$ and all these maps vanish.

\begin{theorem}\label{boundary:exact}
There is a canonical equivariant heart exact sequence
\begin{equation}\label{boundary:heartsequence}
 0\longrightarrow\mathcal N_k\longrightarrow\mathcal M_k[2-k]
 \longrightarrow\mathrm{IC}_k\longrightarrow0
 \qquad(k\in\Z),
\end{equation}
where $\mathcal N_k$ is the complete row $C_l$ with the inherited total shifts, action and differential above.  It is the maximal central subobject, and the quotient is the simple intermediate extension of the prescribed upper incidence object.
\end{theorem}
\begin{proof}
We construct an actual dg morphism, and apply the detector only to bounded coherent objects.  A finite augmented magnetic \v Cech complex applied coefficientwise to $\mathcal M_k$ gives a bicomplex $C^{q,-l}$.  Its ambient actions and matter contractions preserve the vertical degree $q$.  Form the intelligent vertical truncation
\[
 J^{q,-l}=\begin{cases}
 C^{q,-l},&q<2,\\
 \ker(d_v:C^{2,-l}\to C^{3,-l}),&q=2,\\
 0,&q>2.
 \end{cases}
\]
The differential and every homogeneous ambient action preserve these kernels.  With total differential $d_v+(-1)^q d_h$, an ambient element of horizontal degree $t$ acts with the sign $(-1)^{qt}$; thus the inclusion and projection are dg-linear.  Since $H^0_{\{0\}}(T_l)=H^1_{\{0\}}(T_l)=0$, the cycle-to-cohomology map identifies $\operatorname{Tot}(J)$ quasi-isomorphically with the full $H^2$ row placed in vertical degree two.  Both directions are bounded.  This gives the roof
\begin{equation}\label{boundary:roof}
 \mathcal N_k\xleftarrow{\ \sim\ }
 \operatorname{Tot}(J)[2-k]\longrightarrow
 R\Gamma_{\{0\}}(\mathcal M_k)[2-k]
 \longrightarrow\mathcal M_k[2-k].
\end{equation}
The \v Cech terms themselves need not be coherent; the source is the explicit bounded coherent row and the target is bounded coherent.

The construction is equivariant for the finite jet quotient through which the chart data act, together with all loop, source and deck actions.  More precisely use the canonical equivariant functor $R\Gamma_{\{0\}}$ and its counit; a coordinate \v Cech complex is only a computational presentation.  Common refinements preserve the vertical cycles, their projection to $H^2$, and the counit, and identify the roofs in the derived category.  Equivalently they can be formed in an equivariant Cartan--Eilenberg enhancement.  Transport by the action morphism, and the canonical comparisons over its double product, give the algebraic descent cocycle.  Thus no chosen coordinate cover is asserted to be invariant under the jet group, and no averaging over a unipotent group is used.  Naturality for the coefficient maps retains every ambient action and matter contraction in this descent.

The row and $\mathcal M_k$ are termwise flat over the matter polynomial ring.  Their matter-zero restrictions are therefore computed directly; the row differential vanishes there.  After regrading, the row detector consists of the central sheaves $C_l$ in ordinary perverse degree zero.  The detector places $\mathcal N_k$ in the lifted heart, and sends \eqref{boundary:roof} coefficientwise to the canonical map
$H^2_{\{0\}}(T_l)[-2]\to T_l$, shifted by two.
For every finite-length magnetic module $Z$, local-cohomology adjunction and the two vanished lower rows give
\begin{equation}\label{boundary:extidentity}
 \Ext_R^2(Z,T_l)=\Hom_R(Z,H^2_{\{0\}}(T_l)).
\end{equation}
The detected map corresponds to the identity of $C_l$.  It is therefore monic in the ordinary perverse heart, and it contains every central subobject, including nonreduced thickenings.  Regrading weights separate the exterior coefficients.

Exactness and conservativity of the bounded coherent detector now imply that \eqref{boundary:roof} is monic.  For any actual central subobject $Z$ of $\mathcal M_k[2-k]$, the detected composite to its quotient by $\mathcal N_k$ is zero.  An exact conservative functor of hearts is faithful, by applying it to the image of a morphism; the composite itself is therefore zero.  This proves maximality without asserting that arbitrary detected subobjects lift.  A new central subobject of the quotient would contradict this maximality through its preimage.  Proposition~\ref{boundary:depthfailure} excludes all other boundary subobjects and all boundary quotients.  The open restriction is simple, proving the theorem.  At no point is a perverse detector extended to an unbounded or noncoherent local-cohomology object.
\end{proof}

We record the complete equivariant character of this correction.  Put $D=\det V_0$.  For an irreducible target representation $V$ of central weight $w$, let $S(V)$ denote the central simple including its prescribed shift $[-w]$.  External loop/source characters will be written in front.  Then
\begin{equation}\label{boundary:evencharacter}
 \mathcal N_k=u^{2k}S(D^{(k-4)/2})\qquad(k\text{ even}).
\end{equation}
For odd $k$ there is a heart exact sequence
\begin{equation}\label{boundary:oddfiltration}
 0\longrightarrow C_4\longrightarrow\mathcal N_k
 \longrightarrow C_3\oplus C_5\longrightarrow0,
 \quad
 \begin{aligned}
 C_3&=\bigoplus_{j=1}^4u^{2k-1}x_j^{-1}S(D^{(k-3)/2}),\\
 C_4&=u^{2k}S(D^{(k-5)/2}V_0),\\
 C_5&=\bigoplus_{j=1}^4u^{2k-1}x_jS(D^{(k-5)/2}).
 \end{aligned}
\end{equation}
In particular the odd boundary object is generally nonsemisimple; \eqref{boundary:oddfiltration} is a composition filtration, not a direct-sum decomposition.
To verify the weights, set $\delta^2=ab$ and
$W=\operatorname{diag}(a/\delta,b/\delta,1)$.
The point action is $(\mathbf z,\mathbf w)\mapsto(u^{-1}W\mathbf z,u^{-1}\mathbf wW^{-1})$.
The top and bottom rows of $F$ have factors $u^{-1}h/\delta$ and $u^{-3}h/\delta$, whereas the actual quotient rows have $u^{-1}h/x_j$ and $u^{-3}h/x_j$.  Each dual-image generator consequently contributes $x_j/\delta$, with no additional loop factor.  The electric frame contributes $(u^2\delta)^k/x$.
In local cohomology the degree-three and degree-five representatives have loop weight $u^{-1}$; degree four has weight one.  For the latter the source in \eqref{boundary:incidence} has weight $u\cdot u^{-3}$ and the incidence equation supplies $u^2$, so no diagonal shift is missing.  Its even target character is $1$, and its odd target character is $\delta/a+\delta/b$.  Combining the factors proves \eqref{boundary:evencharacter}--\eqref{boundary:oddfiltration}.
The full boundary object satisfies the electric recurrence
$\mathcal N_{k+2}\simeq u^4D\otimes\mathcal N_k$, with the associated central shift $[-2]$: projection onto one imbalance summand identifies either odd eigenspace with the same row matrices, while the retained character changes by $u^4D$.

\subsection{The complete five-label operator}
For the remaining formulas put
\[
 \begin{gathered}
 e_j=e_j(x_1,x_2,x_3,x_4),\quad e_0=1,\ e_4=x,\\
 \delta^2=ab,\quad s^2=a/b,\quad \xi^4=a/b,\quad t=u^2,
 \end{gathered}
\]
\[
 D_d=(1-u^da/b)(1-u^db/a),\quad \mathscr D=D_2D_4,\quad
 h=s+s^{-1},\quad \varrho=s^2+s^{-2}.
\]
Define
\begin{equation}\label{boundary:finiteingredients}
 \begin{gathered}
 c_j=e_j/\delta^j,\qquad
 \mathcal A=\sum_{j=0}^4c_j^2,\qquad
 W_m=\sum_{j=0}^{4-m}c_jc_{j+m},\\
 \mathcal G=\sum_{j=0}^4(-1)^j
       e_j(x_1^2,x_2^2,x_3^2,x_4^2)/\delta^{2j},\qquad
 S_i(u,\xi)=\sum_{\substack{0\leq p,q\leq3\\p-q\equiv i\ (4)}}
                   u^{p+q}\xi^{p-q}.
 \end{gathered}
\end{equation}
The corrected central numerators are the following finite Laurent polynomials:
\begin{equation}\label{boundary:numerators}
\begin{aligned}
 \mathsf E_0&=1+2t+4t^2+2t^3+t^4+t^2\varrho,\\
 \mathsf E_1&=u(1+2t+2t^2+t^3)h,\\
 \mathsf E_2&=t(1+t+t^2)h^2,\\
 \mathsf E_3&=u^3(1+t)h(\varrho+1),\\
 \mathsf E_4&=-1-t^4+2t^2+(t+3t^2+t^3)\varrho+t^2(s^4+s^{-4}),\\[2pt]
 \mathsf O_0&=2t(1+t+t^2)h,\\
 \mathsf O_1&=u(1+t)((1+t)^2+t\varrho),\\
 \mathsf O_2&=th((1+t)^2+t\varrho),\\
 \mathsf O_3&=u^{-1}(1+t)(1+t+t^2)(th^2-1-t^2),\\
 \mathsf O_4&=h(t(1+t)^2\varrho-1-t^4).
\end{aligned}
\end{equation}

\begin{theorem}\label{boundary:operator}
The actual normalized class
\begin{equation}\label{boundary:normalizedclass}
 \mathsf K_k=\frac{u^{2-2k}}{x}[\mathrm{IC}_k],\qquad k\in\Z,
\end{equation}
where $\mathrm{IC}_k$ is defined from $\ell^k/x$ in
\eqref{boundary:heartsequence}, has symmetric coefficients
\begin{equation}\label{boundary:operatorformula}
\begin{aligned}
 f_{2,-2}^{(k)}&=b^{k-8},\\
 f_{1,-1}^{(k)}&=-\frac{u}{xD_4}\sum_{r=0}^4
   e_r b^{k-4-r}\xi^{k-2r}S_{2r-k}(u,\xi),\\
 f_{0,0}^{(k)}&=\frac{u^2\delta^k}{x^2\mathscr D}
 \begin{cases}
 \displaystyle\frac{\mathcal A\mathsf E_0+(-1)^{k/2}\mathcal G D_4}{2}
       +\sum_{m=1}^4W_m\mathsf E_m,&k\text{ even},\\[2pt]
 \displaystyle\frac{\mathcal A\mathsf O_0}{2}
       +\sum_{m=1}^4W_m\mathsf O_m,&k\text{ odd}.
 \end{cases}
\end{aligned}
\end{equation}
The other two coefficients are Weyl transforms.  The formulas descend integrally, are bar-fixed, and satisfy A1--A3.  The middle coefficient is $O(u)$; the central coefficient is $O(u^2)$ for even $k$ and $O(u)$ for odd $k$.  The associated ordinary difference operator preserves the integral symmetric Laurent polynomial ring.  The normalized electric-zero member is $k=8$, equivalently the geometric twist $\ell^8/x^2$ before loop normalization.
\end{theorem}

\begin{proof}
We first derive the characters from the actual image.  For a triple
$L_\varepsilon=(s,s^{-1},\varepsilon)$, $\varepsilon=\pm1$, let $h_n$ be the complete symmetric polynomial and put
\begin{equation}\label{boundary:Hilbert}
 Q_m(L_\varepsilon)=(1-u^2)
      \sum_{n\geq0}u^{2n+m}h_{n+m}(L_\varepsilon)h_n(L_\varepsilon).
\end{equation}
This is the $H^0$ section character of the gauge-charge-$m$ module on the incidence cone: the trace relation removes the bidegree-$(n+m-1,n-1)$ copy from bidegree $(n+m,n)$.  It is not an Euler characteristic in which the intermediate groups have been discarded.  Opposite gauge fibers tensor to the degree-zero trivial fiber, so the section module depends on $|a-b|$, not on $a+b$.
There is a finite rational definition of \eqref{boundary:Hilbert}.  If
$L_\varepsilon=(L_1,L_2,L_3)$ and
$d_i=L_i^2/\prod_{j\ne i}(L_i-L_j)$, partial fractions give
\begin{equation}\label{boundary:partialfractions}
 Q_m(L_\varepsilon)=(1-u^2)u^m
       \sum_{i,j=1}^3\frac{d_id_jL_i^m}{1-u^2L_iL_j}.
\end{equation}
Thus the formulas can be checked by nine rational terms, without expanded high-degree polynomials.

Put $R_m=(Q_m(L_1)+Q_m(L_{-1}))/2$ and
$O_m=(Q_m(L_1)-Q_m(L_{-1}))/2$.
The identity and $\gamma^2$ traces give these even and odd parts.
The traces of $\gamma,\gamma^3$ vanish when the two exterior degrees differ.  When they agree and equal $j$, exchanging the exterior blocks contributes
$(-1)^je_j(x_1^2,\ldots,x_4^2)/\delta^{2j}$.
The neutral-ring trace is $1/D_2$: on bidegree $(n,n)$, the interchange trace is the complete character of
$(s^2,s^{-2},-1)$; the incidence equation has $\gamma$-character $-1$, and hence contributes $1+u^2$.  This cancels the factor $1+u^2$ in its symmetric-algebra denominator.  Consequently the naive sector-$-k$ central character is
\begin{equation}\label{boundary:projector}
 F_k=\begin{cases}
 \mathcal A R_0/2+\sum_{m=1}^4W_mR_m
                  +(-1)^{k/2}\mathcal G/(2D_2),&k\text{ even},\\
 \mathcal A O_0/2+\sum_{m=1}^4W_mO_m,&k\text{ odd}.
 \end{cases}
\end{equation}
This is the actual finite-group projector, including its charged sign.

The upper normalization is independent of lower coefficients.  At depth two,
\[
 \frac{E_{\mathrm{int},2}(B)}{P_2(B)}=\frac{u^{16}x^2}{B^8}.
\]
The intrinsic source $\ell^k/x$ therefore gives raw upper coefficient
\[ u^{16}xB^{k-8}P_2(B)/V_4. \]
Comparing with
\[ J_{2,-2}(A,B)(u^{-2}B)^{k-8} \]
gives precisely the scalar $u^{2-2k}/x$ in \eqref{boundary:normalizedclass}.
At the center there is no further polarization factor; the electric compensation and this scalar give the prefactor $u^2\delta^k/x^2$ in \eqref{boundary:operatorformula}.

The quotient \eqref{boundary:heartsequence} subtracts positive central simple characters.  Indeed $S(V)$ of central weight $w$, including negative $w$, has $K$-theoretic sign $(-1)^w$ from its shift $[-w]$, while simultaneous target negation contributes the same sign.  They cancel.  The analogous shift and exterior signs cancel for the naive object.  Thus no extra alternating sign is attached to the odd exterior degrees in \eqref{boundary:oddfiltration}.
The normalized subtraction is
\begin{equation}\label{boundary:normalizedcorrection}
 \Delta_k=\begin{cases}
 u^2\delta^{k-4}/x,&k\text{ even},\\
 u\delta^{k-3}e_3/x^2+
 u^2\delta^{k-5}(a+b)/x+
 u\delta^{k-5}e_1/x,&k\text{ odd}.
 \end{cases}
\end{equation}
Since $W_4=x/\delta^4$ and
$W_3=e_3/\delta^3+e_1x/\delta^5$, this replaces
$R_4$ by $R_4-1$ in the even projector and replaces
$O_3,O_4$ by $O_3-u^{-1},O_4-h$ in the odd projector.
Applying \eqref{boundary:partialfractions} gives exactly
\eqref{boundary:numerators}: $\mathsf E_m/\mathscr D=R_m$ for $m<4$,
$\mathsf E_4/\mathscr D=R_4-1$, and
$\mathsf O_m/\mathscr D=O_m$ for $m<3$, with the two stated corrections for $m=3,4$.  This proves the central formula by geometry, not by bar averaging or fitting residues.

For clarity we also give the actual middle calculation.  In the $A_3$ cover write
$\alpha=p_0q_0$, $d=p_0^4$, $c=-q_0^4$.
Translating the base but retaining the original matter gives
\[
 E=a_0-\alpha a_1-cb_1,\qquad
 F_2=q_0^2b_0+p_0q_0^3b_1-p_0^2a_1.
\]
The four quotient equations are
$(E,-\alpha E+q_0^2F_2,p_0^2F_2,-p_0q_0^3F_2)$.
Thus the dual image is a free line plus a residual line of deck degree $-2$, not $-1$.
In raw variables set
$\eta^4=A/B$, $P=u^{3/2}\eta$, $Q=u^{1/2}\eta^{-1}$ and
$C=u^{3/2}A^{1/4}B^{3/4}$.
The two dual-generator weights are $u^{-1}x_j/A$ and
$u^3x_j/(CP)$; the electric compensation is $C^k$.
The free-generator Euler product cancels the external constrained mode, and internal polarization inserts $\prod_j(1+ux_j/B)$.
Using $E_{\mathrm{int},1}(B)/P_1(B)=u^4x/B^4$, the common scalar above, and $A=u^{-1}a,B=ub$, yields the middle formula in \eqref{boundary:operatorformula}.  Its congruence is $2r-k$, as required by the actual residual line.

Integral descent in the middle follows because $k-2r+p-q$ is divisible by four.  At the center, $\mathsf E_m$ has $s$-parity $m$, whereas $\mathsf O_m$ has parity $m-1$.  Thus every $\delta^ds^p$ has $d\equiv p\pmod2$ and equals
$a^{(d+p)/2}b^{(d-p)/2}$.  The displayed halves are integral: $\mathsf O_0$ is divisible by two; in the even case substitute
\[
 e_j(x_1^2,\ldots,x_4^2)=
 1,\ e_1^2-2e_2,\ e_2^2-2e_1e_3+2x,\ e_3^2-2e_2x,\ x^2
\]
for $j=0,\ldots,4$.
Modulo two, $\mathsf E_0\equiv D_4\equiv1+t^2\varrho+t^4$, so the half-sum has integral coefficients.
Every numerator in \eqref{boundary:numerators} satisfies
$P(u^{-1},s)=u^{-8}P(u,s)$, while
$\mathscr D(u^{-1},s)=u^{-12}\mathscr D(u,s)$.
The prefactor $u^2$ proves central bar invariance.  The even numerators have nonnegative $u$-degrees.  Only $\mathsf O_3$ has degree $-1$, proving the stated order distinction.  Middle reciprocity follows from
$(p,q)\mapsto(3-q,3-p)$, which preserves the residue condition and complements degree to six.

Finally, \eqref{boundary:normalizedclass} is an actual coherent Coulomb-algebra class.  The residue-image theorem of \cite{SS}, in the conventions of Section~\ref{sec:residues}, therefore proves A3 uniformly for every $k$.  Explicitly, all possible ordinary denominators are distinct simple even-root factors on walls $B/A=u^{2h}$, $-3\leq h\leq3$.  Writing a label as $(R,-R)$, the paired labels are $R$ and $h-R$:
\begin{equation}\label{boundary:walls}
\begin{array}{c|c@{\qquad}c|c}
 h&\{R,R'\}&h&\{R,R'\}\\ \hline
 -3&\{-1,-2\}&1&\{2,-1\},\{1,0\}\\
 -2&\{0,-2\}&2&\{2,0\}\\
 -1&\{1,-2\},\{0,-1\}&3&\{2,1\}\\
 0&\{2,-2\},\{1,-1\}&&
\end{array}
\end{equation}
Self-paired coefficients have no self-wall pole.  For an integral symmetric Laurent polynomial $\phi$, the two shifted values
$\phi(u^{2R}A,u^{-2R}B)$ and
$\phi(u^{2(h-R)}A,u^{-2(h-R)}B)$ agree on their wall, since the arguments are interchanged.  A3 therefore makes the numerator of the full operator value divisible by each wall factor.  Such a factor is primitive and, up to a Laurent unit, monic linear in $B$, so this divisibility holds over the original integer Laurent coefficient ring.  The distinct factors are nonassociate primes in that UFD; their product divides the numerator.  All poles cancel integrally, and Weyl covariance gives a symmetric output.
\end{proof}

The formulas have the exact electric step-four recurrence
\begin{equation}\label{boundary:recurrence}
 f_{R,-R}^{(k+4)}=a^{2-R}b^{2+R}f_{R,-R}^{(k)}
 \qquad(-2\leq R\leq2).
\end{equation}
Indeed, increasing $k$ by four leaves the middle deck index unchanged and
multiplies its prefactor by $b^4\xi^4=ab^3$. At the center it preserves parity
and the even-electric sign $(-1)^{k/2}$, while $\delta^4=a^2b^2$; the top term
gains $b^4$, and Weyl reflection supplies the other labels. The symmetric
multipliers are Laurent units independent of $u$, so bar invariance and lower
order propagate in both electric directions.
This is not literal periodicity and not multiplication of the whole operator by one scalar.  In ordinary variables its multiplier is
$u^{-2R^2}A^{2-R}B^{2+R}$; paired multipliers agree on
$B/A=u^{2h}$ because $R+R'=h$.  Thus direct rational checks of all ten signed residue pairs, self-poles, integrality and bar symmetry in $k=0,1,2,3$ give independent finite checks for every positive or negative integer $k$.  These checks confirm the geometric argument; they do not replace algebra membership or the dg construction.

For $k=8$, the middle coefficient simplifies to
\[
 -\frac{u}{xD_4}\left[
 (1+u^2+u^4+u^6)(a^2b^2+e_2ab+x)
 +(u^2+u^4)(a+b)(e_1ab+e_3)\right],
\]
and the actual central subtraction is $u^2(ab)^2/x$.
For odd $k$ the central expansion instead begins
\[
 f_{0,0}^{(k)}=-u\left(\delta^{k-3}e_3/x^2+
                           \delta^{k-5}e_1/x\right)+O(u^2).
\]
Thus order one is sharp for independent source characters.  Neither the first-strip depth-order estimate nor its positive common-numerator conclusion extends automatically to this mixed support.

\begin{remark}\label{boundary:remaining}
The calculation in this section is specific to four columns and the support
$(2,-2)>(1,-1)>(0,0)$.  On a deeper mixed support, reflexive saturation alone
need not give the heart bound, and lowest local-cohomology defects can require
relative coherence and multistratum gluing.  Section~\ref{sec:limits} gives the
finite coefficient construction, the exact depth thresholds, and a conditional
repair statement with these hypotheses explicit.  In particular an exact
saturated filtration by MCM electric-line pieces cannot be assumed: the depth
lemma would contradict Proposition~\ref{boundary:depthfailure}.  No
all-mixed-sector formula, nonzero first-node-charge theorem, or global
canonical-basis identification follows from the present calculation.
\end{remark}

\section{Mixed-sign sectors and uniqueness}\label{sec:limits}
\subsection{Finite presentations and the exact depth threshold}
Assume four columns and first-node magnetic charge zero. For a general mixed
label $(r,-d)$ with $r>0$ and $d\geq2$, put $\nu=r+d$. On the translated
length-$\nu$ quotient let $T$ be multiplication by $z$. The actual external
matter image is
\begin{equation}\label{limit:image}
 I=\im T^d,\qquad
 \rk I|_{S_j}=\max(\nu-j-d,0)+\max(j-d,0).
\end{equation}
This formula isolates the obstacle: translating the magnetic base does not
translate the regular external-matter subspace into the same incidence problem.

There is nevertheless a uniform finite construction of the naive saturated
object. On a normal affine slice $X=\Spec R$, write $\ell^{[k]}$ for the
reflexive electric extension and use an actual finite matrix for $T^d$. Its
coefficient modules are
\begin{equation}\label{limit:finite-presentation}
 T_l=\left(
 \Lambda^l\bigl((I^*)^{\oplus4}\bigr)\otimes_R\ell^{[k]}
 \right)^{**},\qquad 0\leq l\leq4r,
\end{equation}
with the source and loop factors retained. Both sides of the identification
\[
 T_l=\iota_*\left(
 \Lambda^l\bigl((I_U^*)^{\oplus4}\bigr)\otimes\ell_U^k\right)
\]
are reflexive and agree on the big open regular orbit. This proves the
identification. Matrix images, duals, exterior powers, tensor products, and double
duals are finite-presentation operations. The ambient wedge action and matter
contractions extend from $U$ by Hartogs; their identities hold globally because
the targets are torsion-free. Thus~\eqref{limit:finite-presentation} specifies a
bounded coherent saturated dg object, not only its K-class.

A transverse slice at stratum $j$ has dimension $2j$. Normalize its middle
perversity to be $-j$ on the open stratum and zero at the vertex. Under the
prescribed weight regrading a sheaf coefficient is $T_l[j]$. Derived restriction
lies in degrees at most $-j$, whereas the costalk begins in degree
$\operatorname{depth}_{\mathfrak m}T_l-j$. Hence the two relevant thresholds are
\begin{equation}\label{limit:thresholds}
 \begin{array}{ll}
 \operatorname{depth}_{\mathfrak m}T_l\geq j
      &\text{for the heart inequality},\\
 \operatorname{depth}_{\mathfrak m}T_l\geq j+1
      &\text{for the strict no-subobject inequality}.
 \end{array}
\end{equation}
Maximal Cohen--Macaulay depth $2j$ is sufficient but not necessary. Reflexivity
supplies depth at least two, so it settles the non-strict bound when $j\leq2$,
not when $j\geq3$. Nonzero local cohomology below degree $j$ is an obstruction
to being in the heart; it is not a degree-zero boundary subobject that can simply
be subtracted.

\subsection{A conditional lowest-row repair}
The construction used for $(2,-2)$ has the following precise general form.
It is stated on one transverse vertex; a relative or multistratum application
requires the additional coherence and compatibility indicated afterward.

\begin{proposition}\label{limit:repair}
Let $j\geq1$, and let $\cM$ be a bounded coherent saturated dg object with a bounded coherent,
t-exact conservative zero-matter detector. Suppose that, after weight regrading,
its coefficients at a transverse $2j$-dimensional slice are $T_l[j]$. Assume
\[
 H^q_{\mathfrak m}(T_l)=0\quad(q<j),\qquad
 C_l=H^j_{\mathfrak m}(T_l)\text{ is finite length}
\]
for every $l$. Assume the heart inequalities at all other strata, and that the
detector on objects supported at this vertex is the corresponding central
perverse detector. Then the full row $C_l$, with its induced action and matter
differential, defines a coherent boundary dg object $\cC$ and a canonical heart
exact sequence
\[
 0\longrightarrow\cC\longrightarrow\cM\longrightarrow\cM_{\mathrm{red}}
 \longrightarrow0.
\]
The object $\cC$ is the maximal subobject supported at the vertex, including
nonreduced thickenings, and $\cM_{\mathrm{red}}$ has no subobject or quotient
there. If the other boundary inequalities are strict and the open restriction
is simple, $\cM_{\mathrm{red}}$ is the intermediate extension.
\end{proposition}
\begin{proof}
Take a finite augmented \v Cech model for local cohomology on generators of
$\mathfrak m$, coefficientwise. The wedge action and matter contractions preserve
\v Cech degree. Intelligent vertical truncation through degree $j$ is therefore
a dg sub-bicomplex, quasi-isomorphic to the row $C_l$ in degree $j$ by the assumed
lower vanishing. It gives a dg-linear roof from the row to
$R\Gamma_{\mathfrak m}\cM$ and then to $\cM$. The source is coherent because
the row is finite length, although the individual localized \v Cech terms need
not be coherent. No perverse truncation in an unbounded noncoherent category is
being used.

After the detector and regrading the map is
$H^j_{\mathfrak m}(T_l)[-j]\to T_l$, shifted by $j$. For every finite-length
module $Z$, local-cohomology adjunction and the lower vanishing give
\[
 \Ext_R^j(Z,T_l)=\Hom_R\bigl(Z,H^j_{\mathfrak m}(T_l)\bigr).
\]
The map therefore realizes the entire central subobject, not only its socle.
T-exactness and conservativity make the actual dg map monic and maximal.
A central subobject in the quotient would have central-supported preimage larger
than this maximal one. Quotients of the quotient are quotients of $\cM$, and
the strict derived restriction inequality excludes them. Strictness at the other
strata and simplicity on the open orbit complete the intermediate-extension
argument. The bounded coefficient and \v Cech directions ensure all the stated
roofs are within the relevant bounded constructions.
\end{proof}

For shallow stratifications with $j\leq2$, reflexivity supplies the heart
inequality and codimension-two strata have strict costalk bounds. Under the
local-product hypothesis along the deepest smooth stratum, with coefficients
respecting that product and MCM off the transverse vertex, the lowest degree-two
rows have finite length transversely and glue coherently. The proposition then
gives one-step removal of that deepest-stratum subobject. This is a conditional
construction of the correction, not its evaluated character for every support.

The hypotheses also admit a finite local test. Present $R=S/J$ with $S$ a
polynomial ring in $M$ variables carrying their actual equivariant weights.
A finite $S$-free resolution of $T_l$ and graded local duality give
\begin{equation}\label{limit:duality}
 H^q_{\mathfrak m}(T_l)^\vee
   =\Ext_S^{M-q}(T_l,\omega_S).
\end{equation}
The dual is graded vector-space duality, and $\omega_S$ includes its determinant
character as well as the scalar grading shift. Vanishing and finite length of the
right side test the hypotheses. Induced maps on these Ext modules, followed by
duality, compute the actual row action and differential. Characters alone do not
recover those maps. This is a terminating finite-presentation test for each
explicit input; it is not an assertion that the test always succeeds uniformly
in $r,d$.

\subsection{Why the general mixed formula remains open}
On the open orbit, the cyclic filtration of $Q$ filters $\im T^d$ by line
quotients. Its associated graded can compute an open K-class, but it need not
compute the reflexive extension or its local cohomology. In fact the mixed
example gives a direct obstruction to treating this as an exact saturated
filtration. Each saturated electric-line piece is MCM by
Lemma~\ref{lem:cyclicmcm}. If the saturated exterior coefficients admitted the
corresponding exact finite filtration by such pieces, the depth lemma would
make them MCM. The depth-two degree-four coefficient of
Theorem~\ref{thm:main-boundary} contradicts this conclusion. At least one proposed
saturated sequence must fail to be right exact.

For general $\nu$, several bad strata can interact. An intermediate local-cohomology
row may have positive-dimensional support, and after the first correction the
coefficients need not remain sheaves. Pointwise characters do not prescribe the
relative gluing maps. A full mixed theorem therefore still needs uniform depth
bounds at every stratum together with coherent relative lowest-row complexes and
compatible removal, or a coherent dg-compatible model for the necessary perverse
truncations when lower-than-$j$ defects occur. A geometric resolution theorem
computing those same defects would also suffice. Neither finite presentation nor
the abstract existence of intermediate extension proves these missing statements.
The sharp next input is a uniform calculation of the first nonzero equivariant
Ext module in~\eqref{limit:duality}, including its support and induced maps.

\subsection{Triangular uniqueness on arbitrary rank-two supports}
The restriction to the first strip in Corollary~\ref{cor:unique} can be
removed if one works in the actual integral Coulomb image, rather than
in a residue ansatz. We retain four columns and zero first-node magnetic
charge. All component units are removed in the convention of
\eqref{eq:Jglobal}.

For an integer $v$, set $v_- = \max(-v,0)$. For $R\geq S$ define
\[
 E_{R,S}(A,B)=\prod_{j=1}^4
  \left(\prod_{i=1}^{R_-}(1+u^{2i-1}x_j/A)
  \prod_{i=1}^{S_-}(1+u^{2i-1}x_j/B)\right).
\]
These are the internal constrained-mode factors, after simultaneous target
negation. In particular the formula applies when the two magnetic entries
have different signs. Direct factorization gives
\begin{equation}\label{eq:all-signs-euler-ratio}
 \frac{E_{R,S}(A,B)}{P_{R_-}(A)P_{S_-}(B)}
 =u^{4(R_-^2+S_-^2)}x^{R_-+S_-}
       A^{-4R_-}B^{-4S_-}.
\end{equation}
Indeed each individual ratio
$(1+u^{2i-1}x_j/Z)/(1+u^{-(2i-1)}Z/x_j)$
is $u^{2i-1}x_j/Z$, and
$\sum_{i=1}^d(2i-1)=d^2$.

Here a class supported on an orbit closure means a class represented in
equivariant coherent $G_0$ over that closure, with the prescribed finite-stage
transitions. It does not mean an arbitrary rational difference operator whose
nonzero labels lie in the closure.

\begin{lemma}
\label{lem:actual-maximal-symbol}
Let $\alpha$ be an integral equivariant coherent Coulomb class represented
with finite magnetic support in a fixed determinant component and with
zero first-node magnetic charge. Suppose that $(R,S)$ is its highest
nonzero magnetic stratum. Then its symmetric coefficient $f_{R,S}$,
normalized by \eqref{eq:Jglobal}, is an integral Laurent polynomial in
$u$, the target characters, and the source characters.
\end{lemma}
\begin{proof}
On a fixed magnetic orbit the classical incidence is an equivariant vector
bundle: the dimensions of the permitted matter modes are constant there.
Equivariant coherent homotopy invariance identifies its $G_0$ with the
representation ring of the orbit stabilizer. The unipotent radical does
not change this Grothendieck group; the reductive quotient is the target
torus for $R>S$, and $\GL_2$ for $R=S$, together with the source and rotation
tori. Thus the restriction of $\alpha$ is a finite integral Laurent
stabilizer character $\psi$. At $R=S$ it is Weyl symmetric.

One can first remove any larger strata on which the class restricts to
zero, using the closed--open localization sequence for coherent $G_0$.
The highest ordinary localization coefficient is therefore, up to an
integral Laurent character unit,
\[
 B_{R,S}(A,B)=\frac{E_{R,S}(A,B)\psi(A,B)}{V_{R-S}(A,B)}.
\]
Here the calculation can be made on classical truncation: coherent
$G_0$ of the bounded derived intersection is identified with coherent
$G_0$ of its classical truncation, by its cohomology filtration. On this
single orbit the inclusion of the classical incidence into the internal
matter bundle is an inclusion of vector bundles. Its normal bundle is
exactly the negative-mode quotient. The ordinary Koszul resolution of
this regular linear embedding multiplies every class $\psi$ by its
normal Euler class, namely $E_{R,S}$ after target negation. This argument
is valid for every virtual coherent class by homotopy invariance and
linearity, not only for the structure sheaf. Nonnegative modes contribute
no internal constraint, and fixed-point localization of the smooth
magnetic orbit contributes $V_{R-S}^{-1}$, up to the stated character
unit. Increasing the cutoff adds free matter modes and uses the fixed
transition convention of Proposition~\ref{prop:simple}; it does not add
an excess exterior factor to this classical $G_0$ calculation.
Thus the argument uses neither a boundary-extension formula nor an MCM
assertion about such an extension.
Dividing by $J_{R,S}$ cancels the orbit denominator; equation
\eqref{eq:all-signs-euler-ratio} cancels all matter denominators with an
integral Laurent monomial remaining. Finally the change
$(A,B)=(u^{-R}a,u^{-S}b)$ is a Laurent-ring automorphism. This proves the
assertion integrally, without inverting any integer.
\end{proof}

\begin{theorem}
\label{thm:all-support-uniqueness}
Retain $m=4$ and fix $r\geq s$. Among actual integral coherent Coulomb classes with zero
first-node magnetic charge and support contained in
$\overline{\mathcal O}_{(r,s)}$, there is at most one normalized operator
with a prescribed top symmetric coefficient and Weyl mate, all of whose
strictly lower symmetric coefficients are bar invariant and $O(u)$.
The assertion includes mixed-sign supports, and all integral electric
labels. It does not assert that such an operator exists for every
prescribed top coefficient.
\end{theorem}
\begin{proof}
Subtract two candidates. The top restrictions cancel: restriction to the
top orbit is detected injectively by its stabilizer character, since its
ordinary localization differs from that character by the nonzero factor
$E_{r,s}/V_{r-s}$ and a character unit. Closed--open localization therefore
represents the difference on the magnetic boundary. If it is nonzero,
choose its highest remaining dominant stratum $(R,S)$. There are only
finitely many strata in the specified rank-two orbit closure, ordered
by dominance. Lemma~\ref{lem:actual-maximal-symbol} makes the corresponding
symmetric coefficient a finite Laurent polynomial in $u$ over the
integral Laurent ring of the other characters. The $O(u)$ condition
removes all powers $u^d$ for $d\leq0$, while bar invariance pairs every
positive power with a negative power of the same coefficient. Hence this
Laurent polynomial is zero, contradicting the choice of stratum.
Equivalently one can iterate the localization sequence through the finite
stratification. Injectivity of abelianization then identifies the classes.
\end{proof}

\begin{corollary}\label{cor:mixed-unique}
For every $k\in\mathbb Z$, the five-label operator in
Theorem~\ref{boundary:operator} is the unique actual integral Coulomb
operator on $\overline{\mathcal O}_{(2,-2)}$ with top symmetric coefficient
$b^{k-8}$ whose lower coefficients are bar invariant and $O(u)$.
\end{corollary}
\begin{proof}
Existence and the required lower orders are proved in
Theorem~\ref{boundary:operator}; apply
Theorem~\ref{thm:all-support-uniqueness}. The odd-electric central order
$O(u)$ is sufficient; the stronger even-electric order $O(u^2)$ is not
needed for uniqueness.
\end{proof}

Actual algebra membership is essential in this argument. Residue
conditions alone control root denominators, but do not force the
highest coefficient to contain the negative-mode matter factors
$P_{R_-}(A)P_{S_-}(B)$. Their cancellation above follows from restriction
to the genuine incidence vector bundle. The theorem consequently gives
no converse image theorem, no construction on unresolved mixed supports,
and no exhaustion or multiplicative characterization of a canonical basis.

\begin{remark}\label{rem:matter-divisibility}
Retain four algebraically independent invertible source characters $x_1,\ldots,x_4$.
In determinant component $-4$, consider the nonzero ordinary operator
\[
 \Delta=u^2\,\mathsf U_{-2,-2}.
\]
It is Weyl covariant, has no poles whatsoever, satisfies every paired root-residue
identity, and preserves the integral symmetric Laurent polynomial ring. In the
normalization of \eqref{eq:Jglobal}, its only symmetric coefficient is
\[
 g(a,b)=\frac1{P_2(u^2a)P_2(u^2b)}
 =\prod_{v=1}^4\prod_{Z\in\{a,b\}}
 \frac1{(1+u Z/x_v)(1+u^{-1}Z/x_v)}.
\]
Indeed $J_{-2,-2}(A,B)=u^2P_2(A)P_2(B)$ and
\[ B_{-2,-2}(A,B)=J_{-2,-2}(A,B)g(u^{-2}A,u^{-2}B)=u^2. \]
There is no determinant-component discrepancy: both formulas use the common
component-unit-removed convention. Restoring its common unit multiplies the
ordinary operator and $J$ by the same factor, leaving $g$ unchanged.
Each paired factor of $g$ is fixed by $u\mapsto u^{-1}$. Moreover
\[
 g(a,b)=\frac{u^8x^2}{a^4b^4}
 \prod_{v=1}^4\prod_{Z\in\{a,b\}}
 \frac1{(1+u Z/x_v)(1+u x_v/Z)},\qquad x=\prod_vx_v.
\]
Thus $g$ is bar invariant and has exact positive order $8$, with nonzero leading
coefficient $x^2/(a^4b^4)$, over the full independent-source coefficient ring.
In the closure of the orbit $(-1,-3)$, this is a strictly lower central label.
Therefore $0$ and $\Delta$ have the same prescribed zero top coefficient and
bar-fixed positive-order lower coefficients, disproving residue-only triangular
uniqueness on general supports. The example is outside the actual Coulomb image:
on its highest incidence stratum an actual coefficient must contain the
negative-mode Euler product, while $u^2$ does not. Equivalently the actual-image
maximal-symbol lemma would require $g$ to be Laurent, which it is not.
\end{remark}

\subsection{The domain of triangular uniqueness}\label{subsec:odd-poles}
The two uniqueness results use different hypotheses. Corollary~\ref{cor:unique}
also applies to the restricted even-root residue class on the first strip,
whereas Theorem~\ref{thm:all-support-uniqueness} requires actual algebra
membership on general supports. Even on the first strip, enlarging the
root-pole domain creates a further, distinct ambiguity. The broader
root-binomial coefficient domain in~\cite[Section 3.1]{Gaiotto} permits odd as
well as even ordinary walls, whereas the displayed residue tests inspect only
even walls. For the first-strip support with $N=4$, keep the four noncentral
coefficients fixed and add to the central symmetric coefficient the function
\begin{equation}\label{eq:oddambiguity}
 h(a,b)=\frac{u}{(1-ua/b)(1-ub/a)}.
\end{equation}
This is Weyl symmetric, bar fixed, and $O(u)$. Its ordinary central contribution is
\[
 -u^{-1}h(ua,ub)=-\frac{1}{(1-ua/b)(1-ub/a)}.
\]
It has only simple odd-wall poles
and no residue at any even wall. Thus the broader elementwise tests do not force
this addition to vanish.

More generally, with the noncentral coefficients fixed, the ambiguity in that
broader class is exactly the additive set of Weyl-symmetric bar-fixed $O(u)$
rational Laurent functions with integral Laurent numerator and at most simple
odd-root poles. Necessity follows by subtracting the even-wall equations and
removing the resulting even poles; sufficiency follows because a central shift
preserves $a/b$. Extra electric-degree restrictions, if required, must be stated
separately. These additions need not belong to the actual Coulomb algebra or
preserve symmetric Laurent polynomials. The example is a limitation of a larger
coefficient domain, not a disproof of a complete $K_{\mathfrak q}$-algebra
construction.

The remaining geometric problem is to construct the correction complexes on
general deeper mixed-sign supports. The local Ext test above specifies the
missing data: uniform depth bounds, the supports of the first nonzero defect
modules, and their compatible dg maps. The uniqueness theorem would then
characterize any normalized operators satisfying its hypotheses. The present
results do not address nonzero first-node charge, variable gauge rank, or the
multiplication, trace, and pairing of an entire basis. In particular,
first-strip numerator positivity is not a statement about mixed coefficients
or basis multiplication.

A separate comparison is to transport the four-column classes through the
upper cluster $\mathcal X$-algebra isomorphism of~\cite{SS}, matching its
quantum-torus and coefficient conventions. This would test which classes
are monomials in specified cluster charts and determine their expansions in
a candidate basis. A comparison with cluster $\mathcal A$-monomials requires
an additional compatible $\mathcal A$-structure. Any proposed identification
of the $(2,-2)$ class must reproduce its derived boundary correction for every
electric twist; neither a theta-basis identification nor such a cluster
$\mathcal A$-structure is established here.

\address{School of Mathematics and Statistics, Heilongjiang University,\\
Harbin 150080, P. R. China\\
\email{tangxm@hlju.edu.cn}}

\end{document}